\documentclass[12pt]{article}

\usepackage{amsmath}
\usepackage{amssymb}

\usepackage[margin=0.5in]{geometry}
\usepackage{cleveref}
\usepackage{graphicx}
\usepackage{amsfonts,color,multicol}
\usepackage{subfig}
\usepackage{url}
\usepackage{float}
\newcommand{\norm}[1]{{\Vert #1 \Vert }}

\renewcommand{\Re}{\operatorname{Re}}

\newcommand{\trho}{\rho}

\DeclareMathOperator\erf{erf}

\begin{document}

\title{More efficient time integration for Fourier pseudo-spectral DNS of incompressible turbulence}

\author{David I. Ketcheson\thanks{Authors listed alphabetically.} \thanks{King Abdullah University of Science and Technology (KAUST),
  Extreme Computing Research Center (ECRC),
  Computer Electrical and Mathematical Science and Engineering Division (CEMSE),
  Thuwal, 23955-6900, Saudi Arabia} \and Mikael Mortensen\thanks{Department of Mathematics, University of Oslo Mathematics and Natural Sciences, Oslo, Norway} \and Matteo Parsani\footnotemark[2] \and Nathanael Schilling\thanks{Technische Universität München, Zentrum Mathematik, Boltzmannstr. 3, 85748 Garching, Germany }}

\maketitle

\abstract{
Time integration of Fourier pseudo-spectral DNS is usually performed
using the classical fourth-order accurate Runge--Kutta method, or other
methods of second or third order, with a fixed step size.  We investigate the
use of higher-order Runge--Kutta pairs and automatic step size control based
on local error estimation.  We find that the fifth-order accurate Runge--Kutta
pair of Bogacki \& Shampine gives much greater accuracy at a significantly
reduced computational cost.  Specifically, we demonstrate speedups of 2x-10x for
the same accuracy.  Numerical tests (including the Taylor--Green vortex, Rayleigh--Taylor
instability, and homogeneous isotropic turbulence)
confirm the reliability and efficiency of the method.
We also show that adaptive time stepping provides a significant
computational advantage for some problems (like the development of a Rayleigh--Taylor
instability) without compromising accuracy.
}

\section{Time integration of Fourier pseudo-spectral DNS}
Direct numerical simulation (DNS) is a key tool in improving our
understanding of turbulent flows.  Simulation of turbulence in
the absence of boundaries is essential for understanding
the nature of turbulence itself, and Fourier pseudo-spectral
methods are usually the tool of choice thanks to their high computational
efficiency, scalability, and accuracy.  Because they use basis
functions with global support, for a given number of degrees of freedom (DOFs),
these methods often provide much
greater spatial accuracy
than would be possible with typical finite element, finite volume,
or finite difference methods.

Pseudo-spectral DNS also requires discretization in time.
Explicit time integration is generally preferred because small time
steps are required in order to satisfy the accuracy requirements of DNS.
Most often\footnote{Unfortunately, the selection of a time integrator is considered
such a matter of course that it is not even specified in many
publications.}, the Fourier pseudo-spectral space discretization is
coupled with the well-known fourth-order accurate Runge--Kutta time
discretization \cite{Kutta_rk4_1901} (hereafter RK4)
see for example \cite{yokokawa200216,ishihara2007small}.

While RK4 is a remarkably
useful general purpose integrator, when combined with a spectral
method in space it has the potential to become the main source
of discretization error, unless used with a very small
time step size.  Sometimes lower-order accurate Runge--Kutta or multistep
methods are used; they require even smaller step
sizes in order to provide time accuracy commensurate with the
spectral accuracy obtained in space.
For example, among the available open source codes,
Tarang \cite{verma2012object} includes RK methods of orders one, two, and four;
Turbo \cite{teaca2009energy} uses a 3rd-order RK method, and Philofluid \cite{iovieno2001new}
uses a low-storage fourth-order RK method.
The use of small time step
sizes is especially significant since such simulations are sometimes
run on the largest available supercomputers in order to simulate high Reynolds
number flows in a reasonable wall-clock time.  Switching to a time integrator that
allows larger steps without compromising accuracy or parallel
scalability is a simple change that could yield significant benefits.

More than a century has passed since Kutta's development of RK4
\cite{Kutta_rk4_1901}, and in that time a lot of work has gone into developing
highly accurate and efficient Runge--Kutta methods.  In this work we
explore the application of some of those methods to Fourier
pseudo-spectral DNS, focusing particularly on the fifth-order
method of Bogacki \& Shampine (henceforth BS5) \cite{Bogacki1996}.  We make use
of the open-source code SpectralDNS \cite{spectralDNS}.

We also investigate the usefulness of automatic step size control
based on local error estimation.  This is a well-established technique
for initial value problems.  In the context of Runge--Kutta methods,
local error estimation is performed using a pair of methods that share
a set of common intermediate stages.  The two methods have different
orders of accuracy (herein we consider pairs of fifth/fourth order),
so that the difference between the numerical solutions they provide
serves as an estimate of the error.  We refer the reader to \cite{hairer1993}
for details.

In proposing a change to the time stepping algorithm, care must be taken
to ensure that the accuracy obtained is at least as good as what
would be obtained with RK4.  We show that results from adaptive
fifth-order integrators agree with what is obtained with traditional
methods to very high accuracy, even when considering pointwise differences
of turbulent flow fields.  We also show that adaptive time stepping does
not negatively impact accuracy, even for problems with emergent instabilities.
It is essential to validate these results at reasonably high Reynold's number,
which requires a scalable parallel code like SpectralDNS and a reasonably large
computing resource.  Simulations presented here were run on the
Shaheen Cray XC40 system at KAUST \cite{shaheen_2}.

\section{Pseudo-spectral discretization of incompressible Navier--Stokes}
We consider incompressible fluid flow modeled by the Navier--Stokes equations with a
divergence-free velocity field.  Defining the modified
pressure $P = p + u\cdot u/2$ (where $p$ is the regular pressure), the incompressible
Navier--Stokes equations can be written (see e.g. \cite{pope_book})
\begin{subequations}\label{Incomp-NS}
\begin{gather}
    \nabla \cdot u = 0, \label{divfree}\\
    \frac{\partial u}{\partial t} = u \times \omega - \nabla P + \frac{1}{\mbox{Re}}\nabla^2 u \label{NS1}.
\end{gather}
\end{subequations}
Here $u$ is the velocity, $\omega=\nabla\times u$ is the vorticity, and $\mbox{Re}$ is
the Reynolds number.
%By taking the divergence of both sides of equation \eqref{NS1}, we obtain the
%following equation for the modified pressure: \comment{MP: give reference} \comment{NS:This is how it is done in \cite{mortensen2016}, but there is not much more explanation there either}:
%\begin{gather}
%    \nabla^2 P = \nabla \cdot ( u \times \omega - \mbox{Ri}\rho\mbox{e}_z). \label{pressureEquation}
%\end{gather}

%==============================================================================

\subsection{Spatial Discretisation}
The spatial discretization is based on the traditional Fourier pseudo-spectral
method; the description here follows that in \cite{mortensen2016} and is implemented
in the SpectralDNS package \cite{spectralDNS} in the Python programming language.
We consider uniform Cartesian grids in two or three dimensions.  Let $u$ denote
the approximation to the velocity field on the grid, let $\hat{u}$ denote
its discrete Fourier transform (DFT), and let $k$ denote a wavenumber.
We can eliminate the pressure from \eqref{NS1} by taking the divergence of
both sides.  Then transforming equation \eqref{NS1} to frequency space
gives the following system of ordinary differential equations (ODEs) \cite{mortensen2016}:
\begin{align}
    \partial_t\hat{u}_k & = \widehat{(u \times \omega)}_k - \frac{1}{\mbox{Re}}
    |k|^2 \hat{u}_k - ik\frac{k \cdot \widehat{(u \times \omega)}_k}{|k|^2}. \label{ufdef}
\end{align}
Here $|k|^2 = k \cdot k$. As is usual in the pseudo-spectral approach, the cross
products in equation \eqref{ufdef}
are evaluated in physical space and the result is then
transformed to frequency space (see e.g. \cite{canuto2006basic}).

\subsubsection{Dealiasing}
Dealiasing is done using the 3/2-rule as described in \cite[p.~134]{canuto2006basic}.
For evaluation (in physical space) of the cross products of
equation \eqref{ufdef}, the vectors containing $u$ and $\omega$ are padded with
zeroes so that there are $\frac{3N}{2}$ frequencies in each direction.
Applying the inverse Fourier transform (IDFT) gives variables on a mesh with
$\left(\frac{3N}{2}\right)^3$ points in physical space. After the operation has
been carried out using these variables, the result is transformed back to
frequency space. Then the highest frequencies are truncated to leave $N$
frequencies in each direction.\footnote{Our dealiasing implementation led to
the presence of a nonzero imaginary part of the Nyquist frequency. We do not
expect this to qualitatively affect the character of the solution, since in a
properly resolved simulation the amount of energy in the Nyquist frequency is
very small and the imaginary part is neglected in the inverse FFT.} 

\subsection{Time discretizations}
The main purpose of this work is to compare different explicit Runge--Kutta methods
for the integration of \eqref{ufdef}.  Specifically, we are interested in improvements
that can be achieved by using highly optimized fifth-order pairs.
We compare the following temporal discretizations:
\begin{itemize}
    \item AB2: Second-order Adams-Bashforth.  This is used only in the first test, to give
            an idea of its relative inefficiency compared to higher-order methods.
    \item RK4 \cite{Kutta_rk4_1901}: The classical four-stage fourth-order method.  The method is
        used with a fixed step size.  We also use this method to generate reference
        solutions.
    \item DP5 \cite{Dormand1980}: The well-known 5(4) pair of Dormand \& Prince.
    \item KCL5 \cite[p.~190]{Kennedy2000} Method RK5(4)8[3R+]M from the work of Kennedy, Carpenter, \& Lewis.
    \item BS5 \cite{Bogacki1996}: The 5(4) pair due to Bogacki \& Shampine.  We
                            use the error estimator $\hat{b}$ (not $b^*$).
\end{itemize}
For the embedded pairs, we compare implementations with fixed step size and
with variable step size based on local error control.

\subsubsection{Automatic step-size control}
We use an automatic step-size control method
based on \cite[p167]{hairer1993}, the details of which we
describe below for the three-dimensional case; the
two-dimensional case is almost exactly the same.

Two different Runge--Kutta methods are used for each step, the second being an
``embedded'' method for automatic step-size control. Let the numerical solution
(in frequency space, with a real DFT in one direction) at timestep $t \in
\mathbb{N}$ be denoted by $\hat{u}^{(n)}$ for the main Runge--Kutta method and
by $\hat{v}^{(n)}$ for the embedded method. As both $\hat{u}^{(n)}$ and
$\hat{v}^{(n)}$ are in frequency space, they each have three components at
every point.  We write $\pi_j : \mathbb{R}^3 \rightarrow \mathbb{R}$ for
$j=1,2,3$ as the canonical projections onto the components (i.e. $u^{(n)} =
(\pi_1u^{(n)}, \pi_2u^{(n)}, \pi_3u^{(n)})^T$). For each wavenumber $k$
considered, the vector $\widehat{sc}_k^{(n)} = (\pi_1 \widehat{sc}_k^{(n)},
\pi_2 \widehat{sc}_k^{(n)}, \pi_3 \widehat{sc}_k^{(n)})$ is defined to have
dimensions like that of $\hat{u}$ and that
\begin{align*}
    \pi_j \widehat{sc}_k^{(n)} = TOL_{abs} + \max\left\{|\pi_j\hat{u}^{(n-1)}_k|, |\pi_j\hat{u}_k^{(n)}|\right\}TOL_{rel},
\end{align*}
where $TOL_{abs}$ and $TOL_{rel}$ are absolute and relative tolerances
respectively. For all the plots in this document, we used $TOL_{abs} =
TOL_{rel}$ with values in the range $10^{-10}$ to $10^{-2}$.
The error for each component $j$ is estimated by
\begin{align*}
    err_j^{(n)} &=\sqrt{\frac{1}{N^2(N/2 + 1)}\sum_k\left(\frac{\pi_j u^{(n)}_k - \pi_j v^{(n)}_k}{\pi_j \widehat{sc}_k^{(n)}}\right)^2}.
\end{align*}
Here the $N^2(N/2 + 1)$ term gives the total number of frequencies $k$
considered. Note that due to the use of the real DFT in the first direction,
this treats the zero and Nyquist-frequencies slightly differently than if a
complex DFT had been used in each direction. The final error estimate is then
formed by taking the maximum over all components, i.e. $err^{(n)} =
\max\{err_1^{(n)},err_2^{(n)}, err_3^{(n)}\}$

Let $h^{(n)}$ be the stepsize used for the $n$-th step, with $h^{(0)}$ being
supplied at the beginning.  We set
$$h_\textup{new} = h^{(n)}\cdot \min\left\{\delta_\textup{max}, \max\left\{\delta_\textup{min}, \delta\cdot (1/err^{(n)})^{1/(q+1)}\right\}\right\},$$
with
$\delta_\textup{max}=2$ and $\delta_\textup{min}=0.01$, $\delta=0.8$ and $q$ being the order of the embedded
method. If $err\leq 1$ then we continue to the next step with $h^{(n+1)} =
h_{new}$. Else the current step is rejected and we re-run the step with
$h^{(n)} = h_{new}$. We set $\delta_\textup{max}=1$ for the first step after a rejected
step.

\section{Numerical comparisons}
We test each of the candidate integrators on variations of three classical problems:
the Taylor--Green vortex, the Rayleigh--Taylor instability, and homogeneous isotropic
turbulence. We calculate (with RK4 and a very small timestep) a reference solution against
which other solutions are compared. As a measure for the error, we use the discrete
$L_2$ and maximum norms of the difference (in velocity or density) between a given solution
and the reference solution.  This is a very stringent test, given that in each case
we are computing a (chaotic) turbulent flow field.  

It is natural to wonder if the use of much larger time steps (allowed by higher-order
methods and by time step adaptivity) might somehow damp out fine features of the flow.
In the worst case, there might be a feedback effect in which numerical dissipation from
the use of large step size with an adaptive integrator prevents the development
of an instability, which in turn causes the integrator to continue with a large step size.
We will see that this problem does not arise in practice.
The most important point to take away
from the results that follow is that the use of high-order, adaptive methods yields essentially
the same solution but requires a much smaller number of FFTs.

\subsection{Taylor--Green vortex}
The Taylor--Green (TG) vortex test case is a widely-used benchmark, see for
instance \cite{de_wiart_2014}.
We solve the system of equations \eqref{Incomp-NS} on the periodic cube $[-\pi \le x,y,z \le +\pi]$.
The initial velocity components are
\begin{align*}
    u^{(0)} = \left(\begin{array}{c}
    \sin{\left(x\right)}\cos{\left(y\right)}\cos{\left(z\right)} \\
    -\cos{\left(x\right)}\sin{\left(y\right)}\cos{\left(z\right)}\\
	0
\end{array}\right).
\end{align*}
The Reynolds number for this flow is
defined as $Re = 1/\nu$, where $\nu$ is the dynamic viscosity. Starting from the initial condition,
the nonlinear interactions of different flow scales yield vortex breakdowns.
This nonlinear process is initially laminar, but it subsequently
develops into near  anisotropic turbulence that decays with the typical
spectral energy distribution.
We consider three Reynolds numbers:
$Re=280$, $Re = 800$ and $Re = 1,600$.

To measure accuracy, we look at the maximum (over all time-steps) absolute
difference (between a computed solution and the reference solutions computed on a fine grid) in the rate
of dissipation of kinetic energy.
This maximum is obtained by evaluating a piecewise quadratic interpolation
of the reference rate of dissipation of kinetic energy at each timestep. The interpolation is done
using the scipy package \cite{scipy}.\\
We obtained qualitatively very similar results by considering the $L^2$-error of the velocity field at the final time.
To calculate the rate of dissipation in kinetic energy, we first note that on the cube $[-\pi,\pi)^3$ we have that
the total kinetic energy per unit volume is given by $E_{kin} = \frac{1}{2N^6}\norm{\hat{u}}$
where $\norm{.}_2$ is the $L_2$ norm, i.e. $\norm{\hat{u}}_2^2 = \sum_k\sum_{j = 1,2,3}|\pi_j \hat{u}_k|^2$
where $k$ ranges over all (discrete) frequencies considered and $\pi_j$ is the canonical projection onto the $j$-th component so that $\hat{u}_k = (\pi_1\hat{u}_k, \pi_2\hat{u}_k, \pi_3\hat{u}_k)^T$.
The rate of dissipation of kinetic energy is given by
$$\epsilon = -\frac{d}{dt}E_{kin} = \frac{-1}{2N^6}\frac{d}{dt}\sum_{k}\sum_{j=1,2,3}(\pi_j\hat{u}_k)(\overline{\pi_j\hat{u}_k}) = \frac{-1}{N^6}\sum_k\sum_{j=1,2,3}\Re\{\pi_j(\hat{u}_kf(\hat{u}_k)) \},$$
where we use the fact that $\frac{d\hat{u}}{dt} = f(\hat{u})$ with the sum ranging over every frequency $k$ considered. Since
we already calculate an approximation of $f(\hat{u})$ when applying the time-integration method, this gives a cheap way of calculating the rate of
dissipation of kinetic energy in frequency space.
%Our results are in excellent agreement with the rate
%of dissipation from the reference \cite{de_wiart_2014}.
Note that in the
implementation we used, the method described above had to be adjusted to treat the 0 and Nyquist frequencies
for the first real DFT correctly.

\subsubsection{Re=280}
\begin{figure}
      \centering
      %\subfigure{\includegraphics[width=0.8\textwidth]{Figures/combined_plots_280/udiff_vs_f_evaluations_edited.pdf}}
      \includegraphics[width=0.8\textwidth]{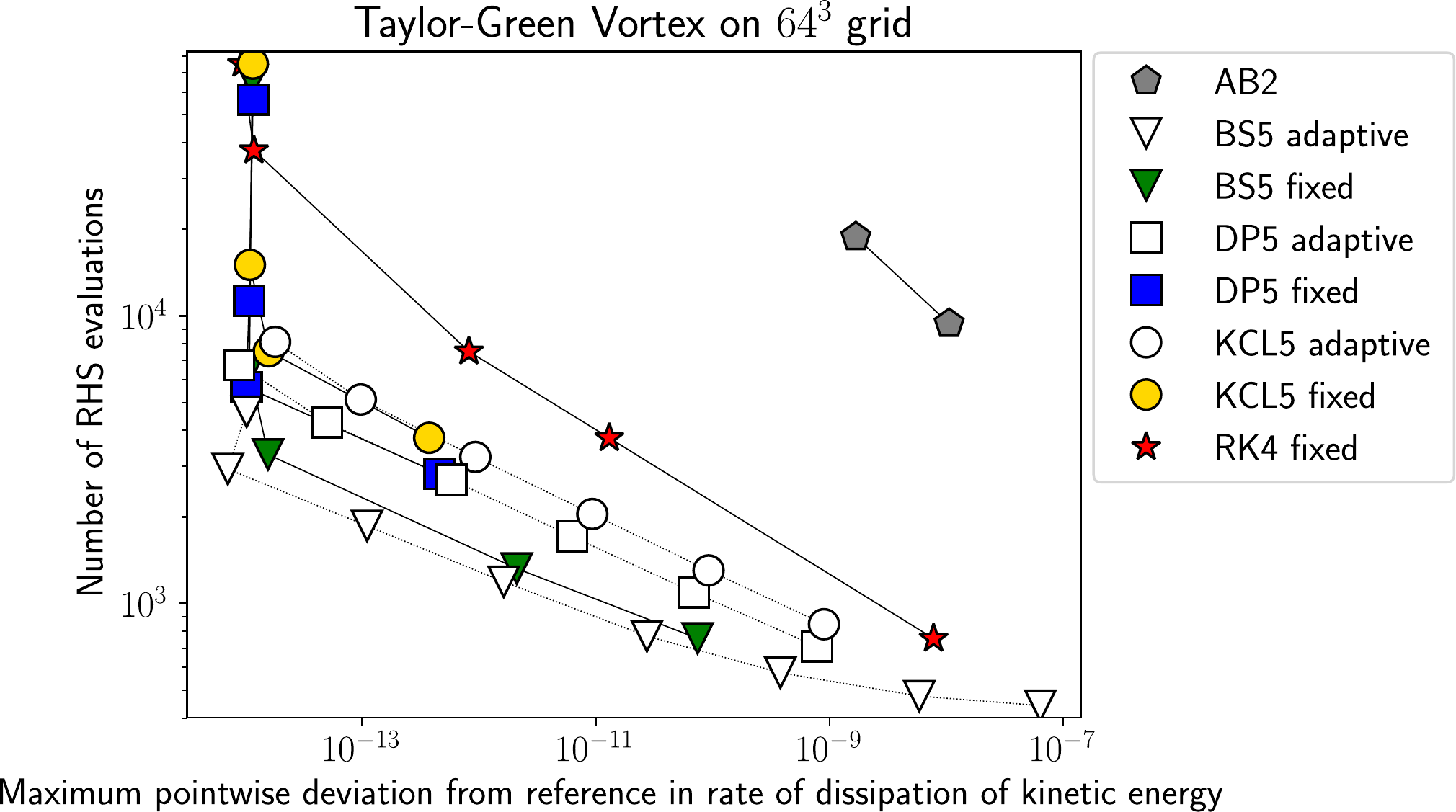}\vspace{1em}\\
      %\caption{Comparison of time integration methods at $\mbox{Re}=280$.}\label{TG280}
      %\subfigure{\includegraphics[width=0.8\textwidth]{Figures/combined_plots_280/udiff_vs_f_evaluations_edited.pdf}}
      \centering
      \caption{Comparison of time integration methods at $\mbox{Re}=280$.}\label{TG280}
\end{figure}

Reference values for the rate of dissipation of kinetic energy are available
from \cite{chapelier2012final}. Here a $64^3$ grid was used, which is sufficient for
DNS \cite{chapelier2012final}.
\Cref{TG280} shows the error in the calculated rate of dissipation of kinetic energy until $T = 9.4$.
%and the maximum norm (on gridpoints) of the error of the solution in physical space at time $T = 9.4$.
The comparison is against a solution generated with
RK4 and a timestep of $10^{-4}$. A timestep of $10^{-3}$ was used for the Re=1600
case on a $512^3$ grid by \cite{mortensen2016}, and hence we expect a timestep of
$10^{-4}$ to be more than adequate for computing the reference solution.

We used a range of tolerances from $10^{-10}$ to $10^{-3}$.
It is clear that that BS5 outperforms the other methods by a wide margin. In particular,
the BS5 with adaptive timestep requires 2-10 times less RHS evaluations than the
RK4 method for a wide range of error norms.

\subsubsection{Re=800}
\begin{figure}
\centering
{\includegraphics[width=0.8\textwidth]{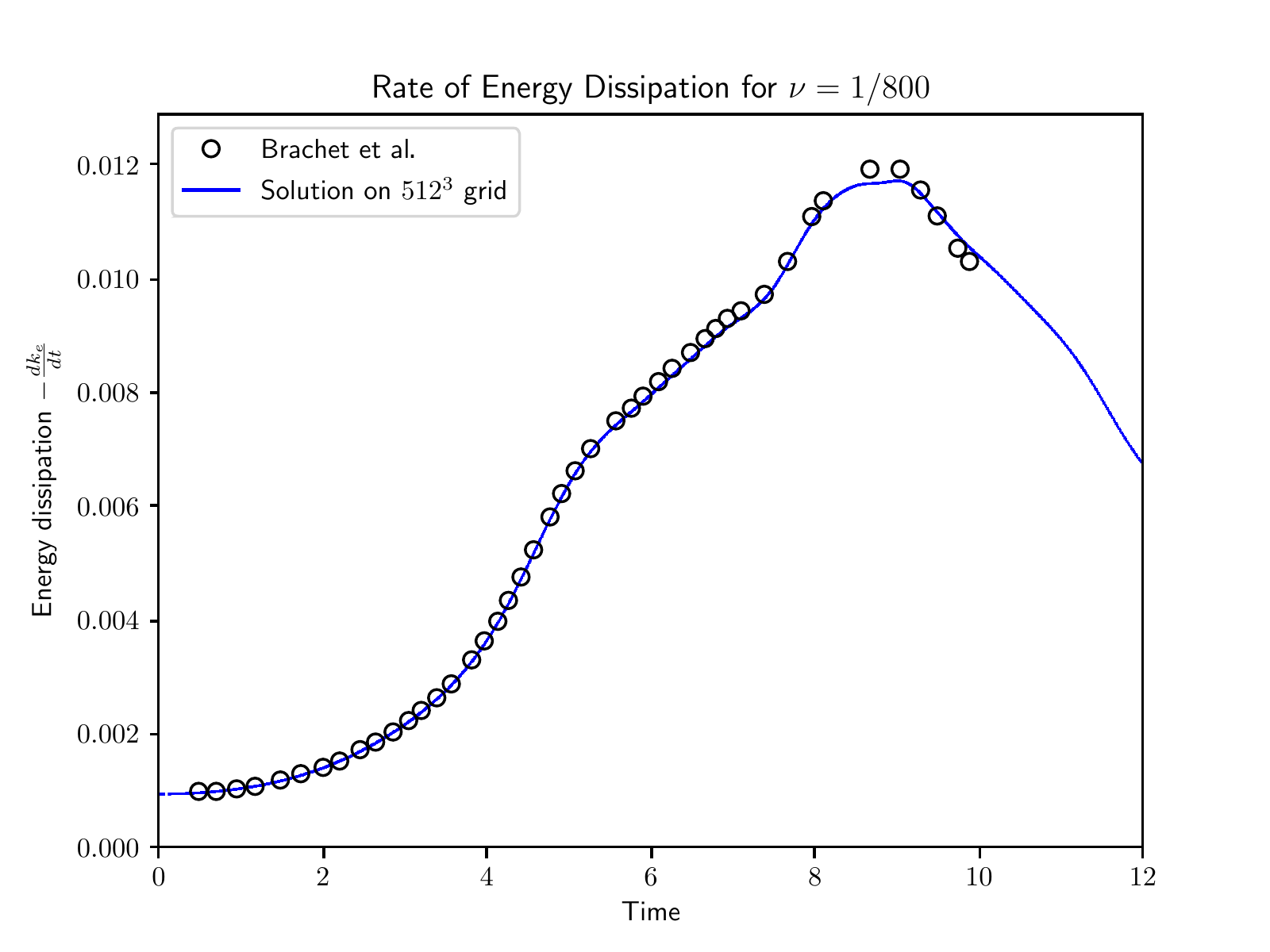}}
\caption{Time rate of change of energy dissipation at $\mbox{Re}=800$: comparison with results of Brachet et al. \cite{brachet1983small}.}
\label{fig:brachetcomparison}
\end{figure}

We used a $256^3$ grid for the Re=800 runs, as done by \cite{brachet1983small}
which is used as a reference by \cite{gassner2013accuracy}.
A comparison of our results
with those reported in \cite{brachet1983small} using a grid with $256^3$ elements
is plotted in Figure \ref{fig:brachetcomparison}.
Our solution 
did not match exactly with that reported by \cite{brachet1983small}. However, the
discrepancy did not go away even when running on a $512^3$ grid.  Similar discrepancies
appear in \cite{gassner2013accuracy}; we refer to that work for a more detailed
discussion.

%\comment{NS:
%I remember the discrepancies with Brachet were a lot smaller
%when using 2/3-rule dealiasing, but 2/3-rule dealising isn't working at the moment so I can't double check this.
%Also, it seems like Brachet doesn't explicitly say he is doing DNS, but only
%that using a $256^3$ grid is accurate for $Re \leq 3000$ \emph{except for the
%smallest scales}.}
%\comment{MP: I think $256^3$ is a good resolution. However, the small scales are
%underresolved.}

\begin{figure}
      \centering
      %\subfigure{\includegraphics[width=0.8\textwidth]{Figures/combined_plots_800/udiff_vs_f_evaluations_edited.pdf}}
      \includegraphics[width=0.8\textwidth]{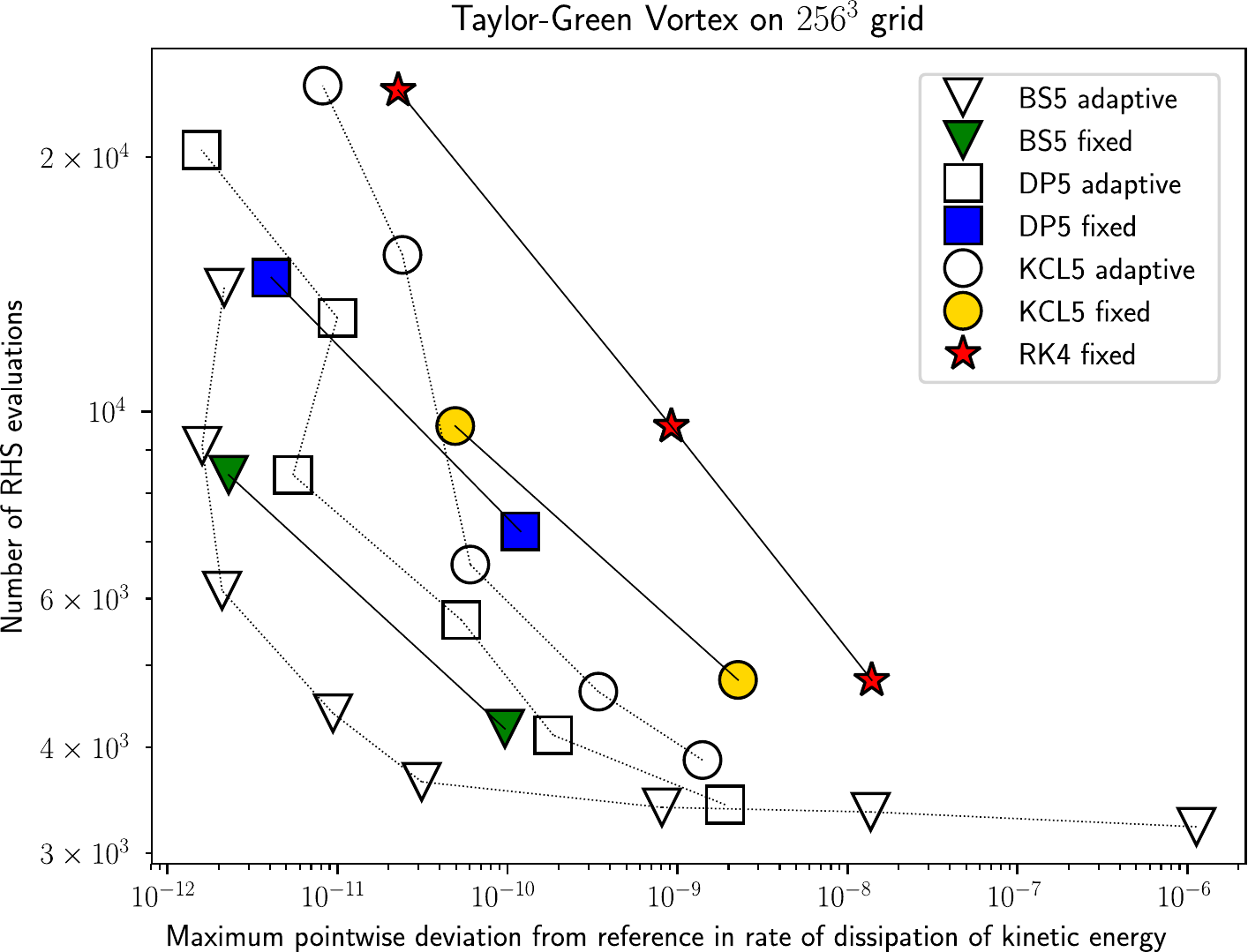}
      %\centering
      %\includegraphics[width=0.8\textwidth]{figs/TG_inferror_800.pdf}
      \caption{Comparison of time integration methods at $\mbox{Re}=800$}\label{TG800}
\end{figure}

We used a reference timestep of $10^{-3}$ and a range of tolerances from $10^{-9}$ to $10^{-2}$. Results are shown in Figure \ref{TG800}.  All 5th-order methods outperform RK4, and the best results are obtained with the adaptive BS5 method, as was found for the TG vortex at $\mbox{Re}=280$.

\subsubsection{Re=1,600}
We used a $512^3$ grid, and a reference timestep of $10^{-3}$ as done in
\cite{mortensen2016}. The TG vortex at this Reynolds number is a well-known test-case, and reference data are
available for the the rate of dissipation of kinetic energy, \cite{de_wiart_2014}.

\begin{figure}
      \centering
      %\subfigure{\includegraphics[width=0.8\textwidth]{Figures/combined_plots_1600/udiff_vs_f_evaluations_edited.pdf}}
      \includegraphics[width=0.8\textwidth]{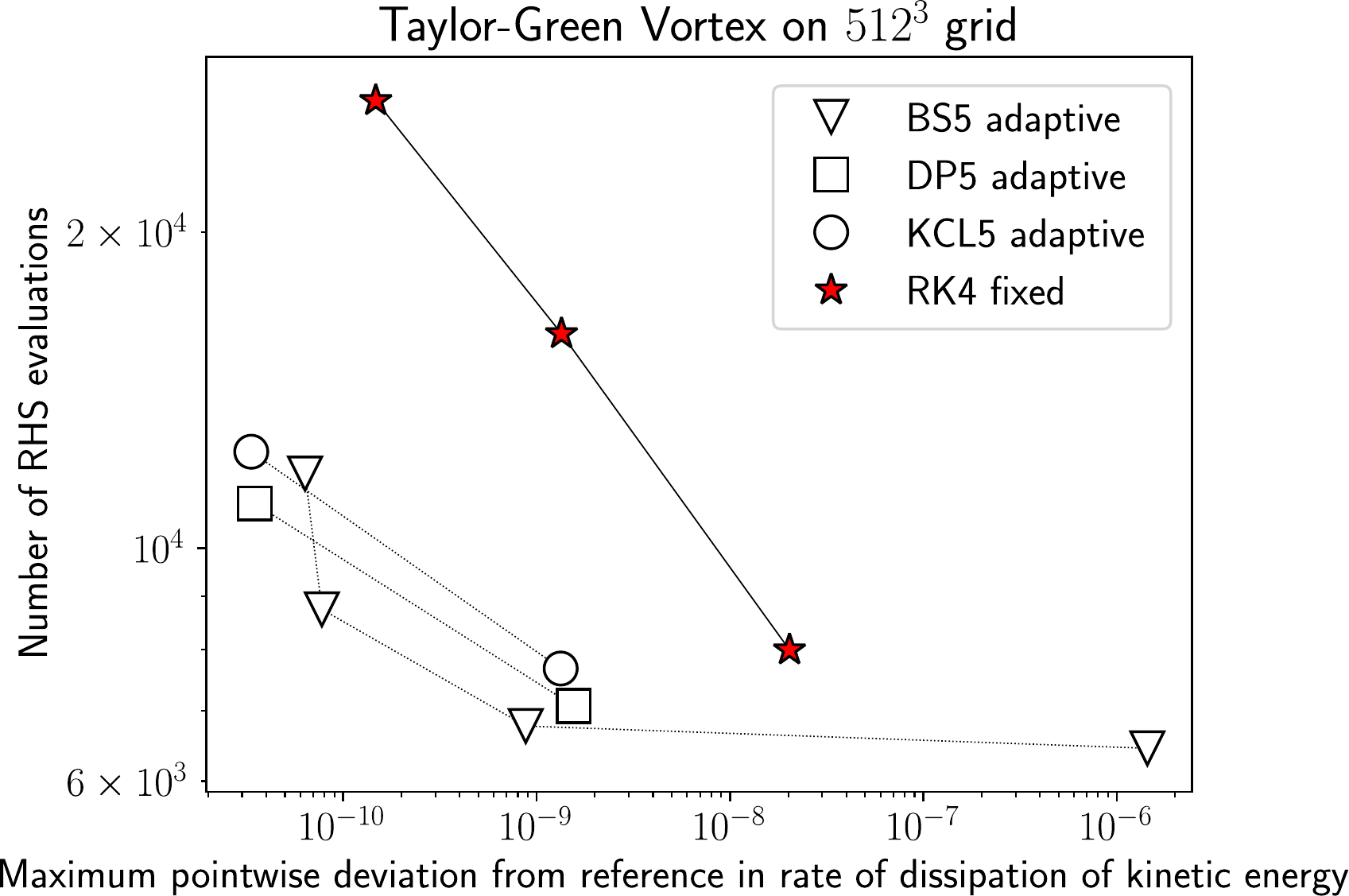}\vspace{1em}\\
      \caption{Comparison of time integration methods at $\mbox{Re}=1600$}\label{TG1600}
\end{figure}

We used a range of tolerances between $10^{-7}$ and $10^{-2}$. Results are shown in Figure \ref{TG1600}.  All 5th-order methods outperform RK4
and, once more, the best results are obtained with the adaptive BS5 method.

\subsection{Rayleigh--Taylor instability}
In this section we solve the incompressible Navier--Stokes
equations in the presence of gravity, with variable density, using the
Boussinesq approximation:
\begin{subequations}
\begin{gather}
    \nabla \cdot u = 0, \\
    \frac{\partial u}{\partial t} = u \times \omega - \nabla P + \frac{1}{\mbox{Re}}\nabla^2 u -
    \mbox{Ri}\trho\mbox{e}_z  \label{NSb1},  \\
    \frac{\partial \trho}{\partial t} = -\nabla\cdot(\trho u) + \frac{1}{\mbox{Re}\mbox{Pr}}\nabla^2
    \trho. \label{NSb2}
\end{gather}
\end{subequations}
Here $\trho$ % = \rho(x,y,z,t)-\overline{\rho}$
denotes the deviation from ambient density,
Ri and Pr denote the Richardson number and the Prandtl number, and $\mbox{e}_z$
indicates the unit vector in the $z$-direction.  In place of \eqref{ufdef},
the system of ODEs in this case takes the form
\begin{subequations}
\begin{align}
    \partial_t\hat{u}_k & = \widehat{(u \times \omega)}_k - \frac{1}{\mbox{Re}}
    |k|^2 \hat{u}_k - ik\frac{k \cdot \widehat{(u \times \omega)}_k -
    \mbox{Ri}\hat{\rho}_k e_z}{|k|^2} - \mbox{Ri}\hat{\rho}_k\mbox{e}_z, \label{ufdefboussinesq}\\
    \partial_t\hat{\rho}_k & = -i k \cdot \widehat{(\rho u)}_k -
    \frac{|k|^2\hat{\rho}_k}{\mbox{Re}\mbox{Pr}}. \label{rhofdefboussinesq}
\end{align}
\end{subequations}

The Rayleigh--Taylor instability is a classical fluid-dynamical instability
that appears in the presence of gravity when a heavier fluid lies above a lighter
fluid.  It is of interest here for two reasons.  First, it is a widely-simulated
phenomenon and serves as a benchmark.  Second, it presents an opportunity to
greatly improve time-stepping efficiency through automatic step size control.
This is because the instability develops slowly at first, so it is expected that
very large step sizes can be used initially.

We simulated the single-mode Rayleigh-taylor instability in both 2D and 3D. The
initial conditions are based on those of \cite{livescu2011direct}, with a slightly smoothed
interface between the fluids in order to allow it to be represented accurately in the
Fourier basis:
\begin{equation}
\begin{aligned}
    u_0 &= 0, \\
    \rho_0 &= \frac{1}{2}(\erf(z - z_0 + \zeta(x,y)))\Delta\rho.
\end{aligned}
\end{equation}
Here $u_0$ and $\rho_0$ are the initial velocity and density fields
respectively and $\erf$ is the error-function; $z_0$ is the location of the
fluid interface.
The $\zeta(x,y)$-term represents a small perturbation to the interface, designed
to seed a single-mode instability.  The quantity $\Delta \rho$
represents the difference in density between the two fluids; we take $\Delta \rho = 1/10.$
This yields an Atwood number of less than $0.05$ for which
the Boussinesq approximation is reasonable. Both the
Richardson number (Ri) and the Prandtl number (Pr) are set to $1$.

\subsubsection{2D}
\begin{figure}
\centering
\subfloat[$T = 0.0$]{\includegraphics[width=0.2\textwidth]{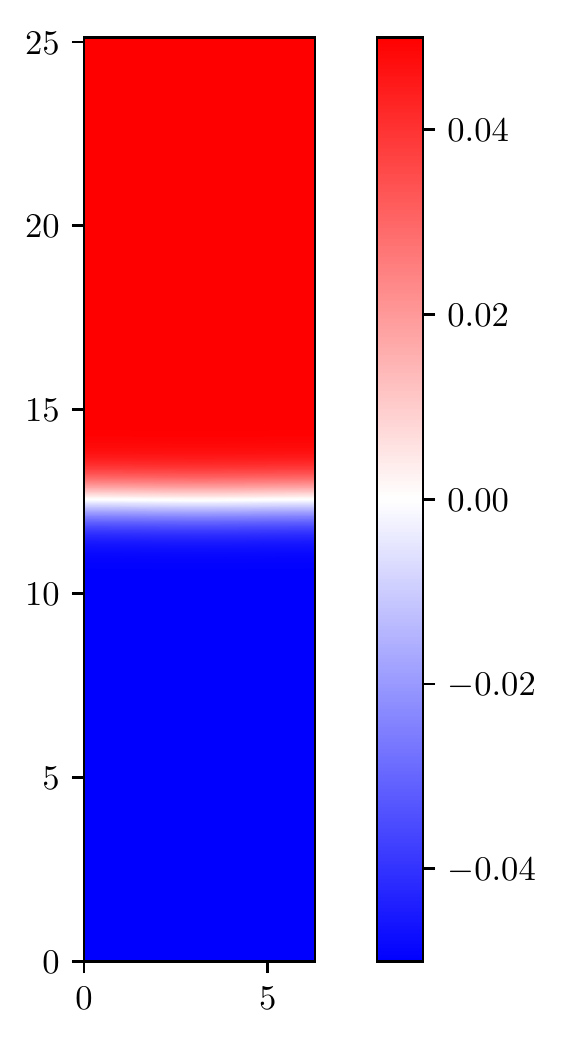}}
\subfloat[$T = 30.0$]{\includegraphics[width=0.2\textwidth]{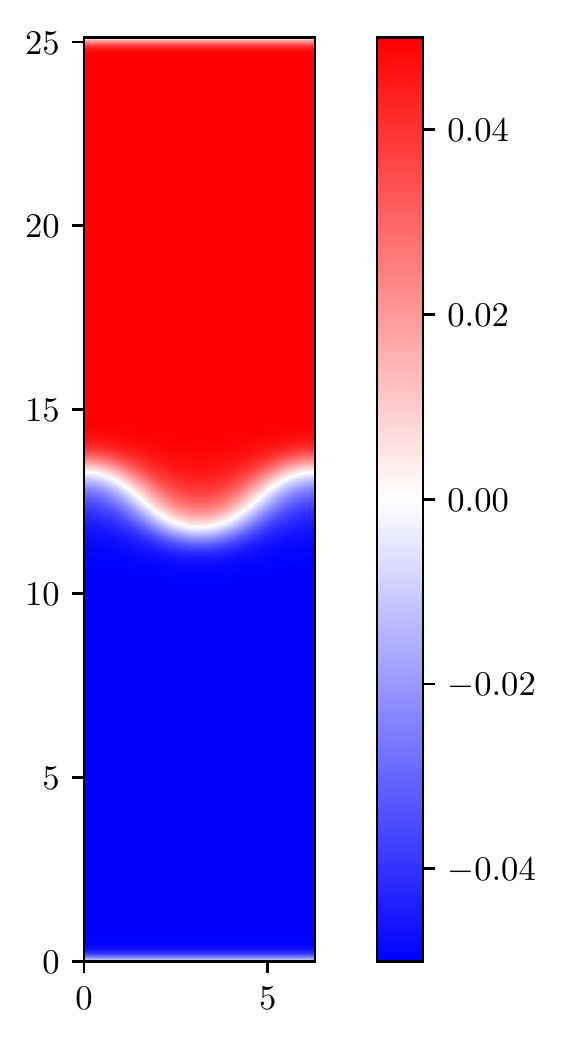}}
\subfloat[$T = 50.0$]{\includegraphics[width=0.2\textwidth]{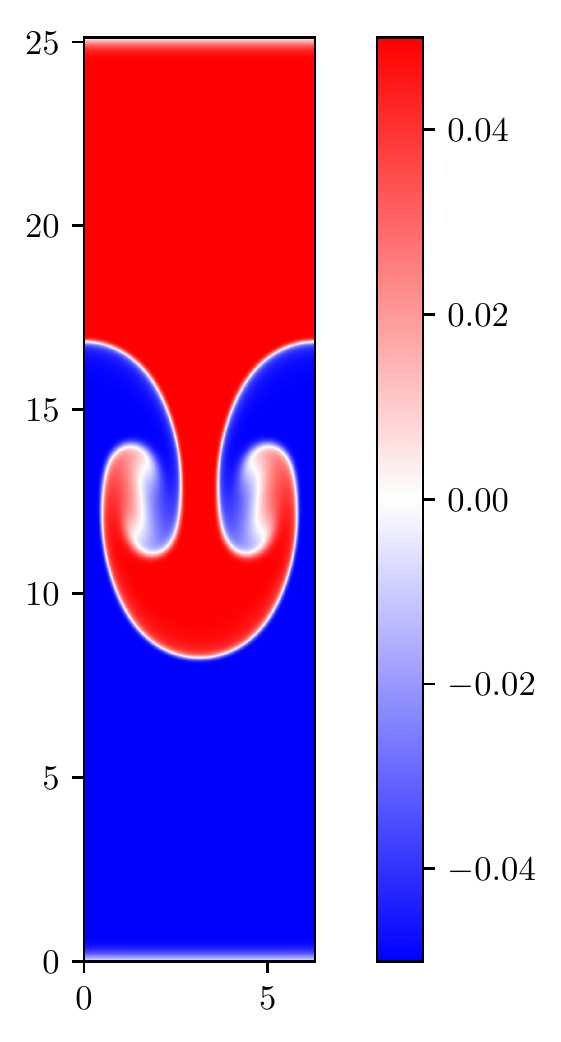}}
\subfloat[$T = 70.0$]{\includegraphics[width=0.2\textwidth]{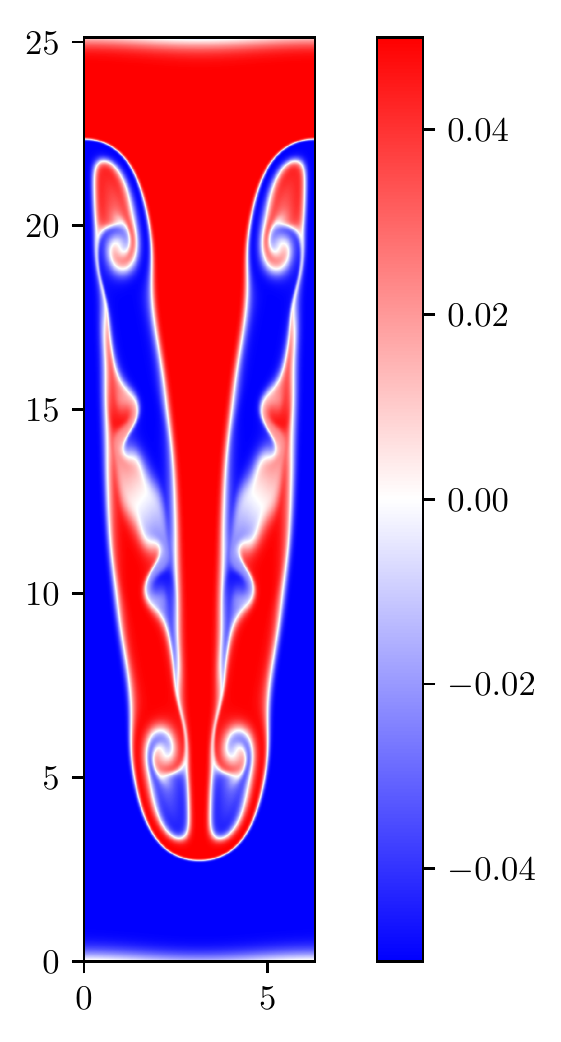}}
\caption{Density field of Rayleigh--Taylor instability on $512\times2048$ grid}
\label{fig:rti2ddensity}
\end{figure}

We use a $512\times2048$ grid with $\nu =  1/1600$
and take $\zeta=-0.01\cos(x)$.
For the reference run, we used a fixed timestep of $10^{-4}$.
To verify that the
requirements for DNS are met, the simulation was also run on a
$1024\times4096$ grid (using BS5 with adaptive timestepping with a tolerance of $10^{-5}$)
and the results compared to the reference on the coarser grid.
The relative $L_2$ and maximum norm of the difference in the velocity field
were both less than $1\%$. For the density, the (relative) $L_2$ norm of the difference
was less than $1\%$ and in the maximum norm the difference was less than
$10\%$. Comparing the kinetic energy and the rate of dissipation (here
negative) of kinetic energy on both grids gives an error of the order $10^{-6}$
and $10^{-10}$, respectively. Differences in the bubble-height between the two
runs are smaller than the resolution of the grid.
The simulation was run from $T=0$ to $T=70.0$. The density field is plotted for various times in Figure \ref{fig:rti2ddensity}.

\begin{figure}
      \centering
      %\subfigure{\includegraphics[width=0.8\textwidth]{Figures/combined_plots_rt_2d/udiff_vs_f_evaluations_edited.pdf}}
      \includegraphics[width=0.5\textwidth]{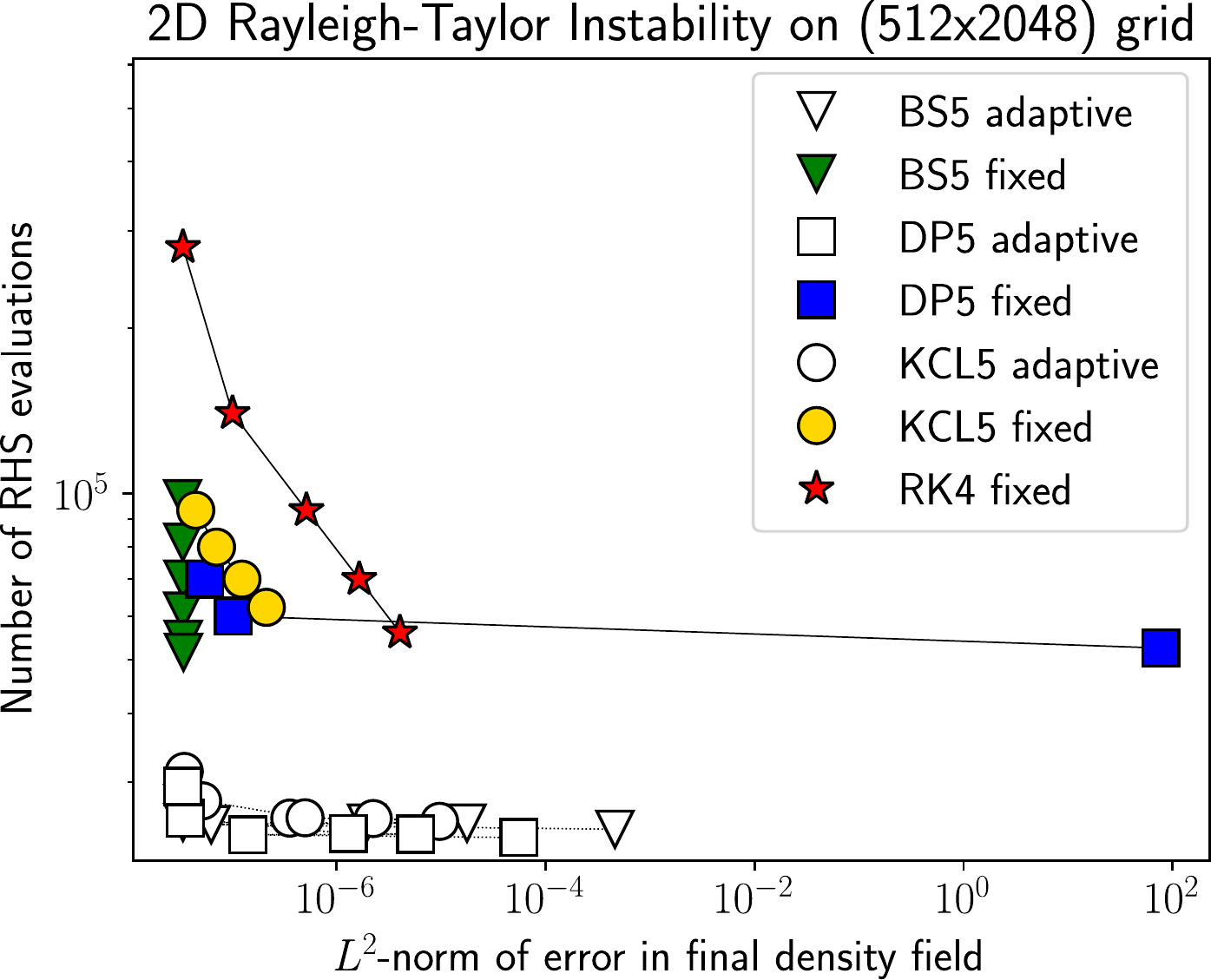}
      \includegraphics[width=0.5\textwidth]{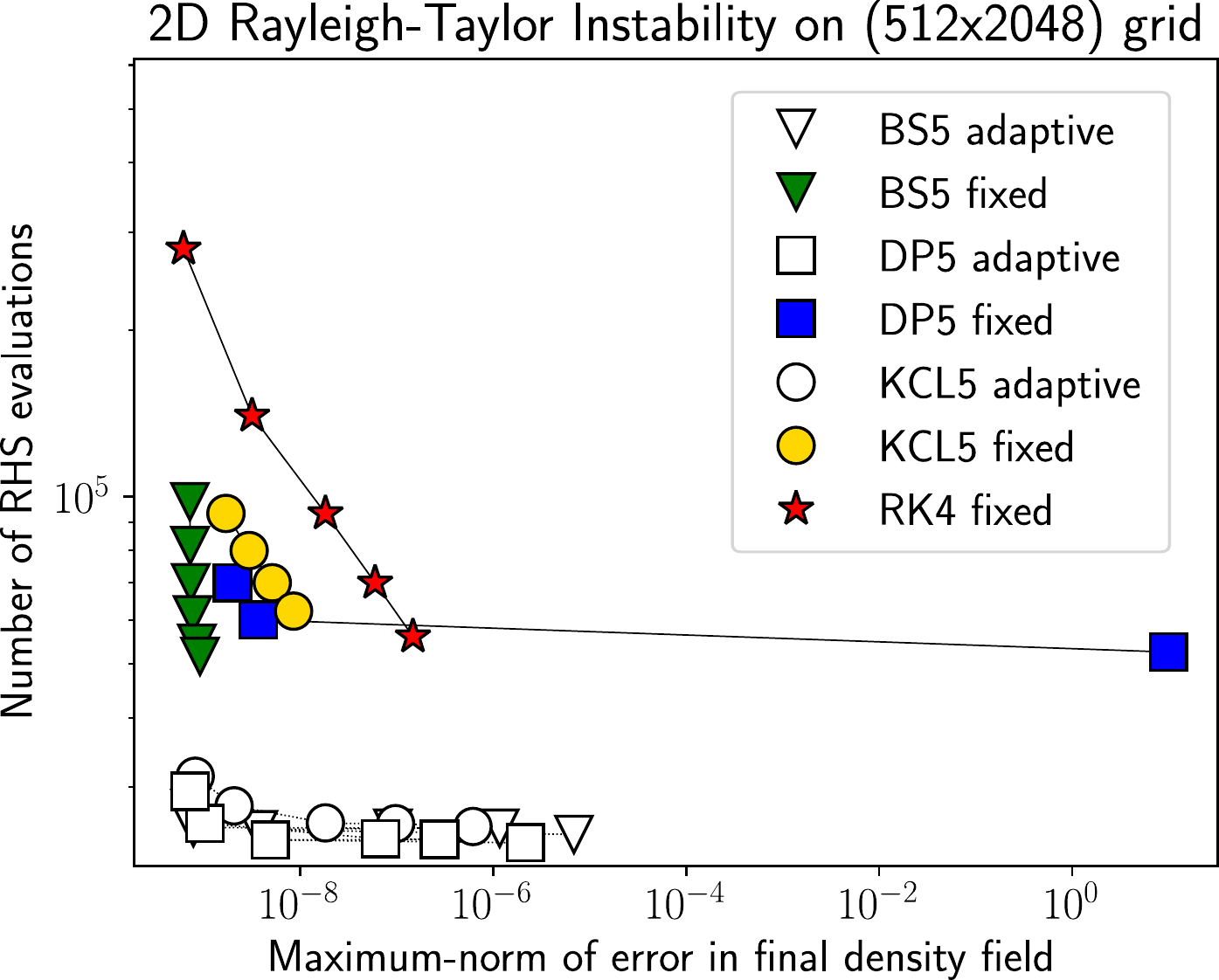}
      \caption{Comparison of time integration methods at $\nu^{-1}=1600$}\label{RTI2D}
\end{figure}

We used a range of tolerances from $10^{-9}$ to $10^{-4}$. Results for all methods are shown in Figure \ref{RTI2D}.
Clearly all the schemes with adaptive timestepping outperform the schemes with fixed timesteps.
Among the
methods with adaptive timestep, the DP5 schemes is the most efficient for this problem.

\subsubsection{3D}
\begin{figure}
\centering
\subfloat[$T = 0.0$]{\includegraphics[width=0.2\textwidth]{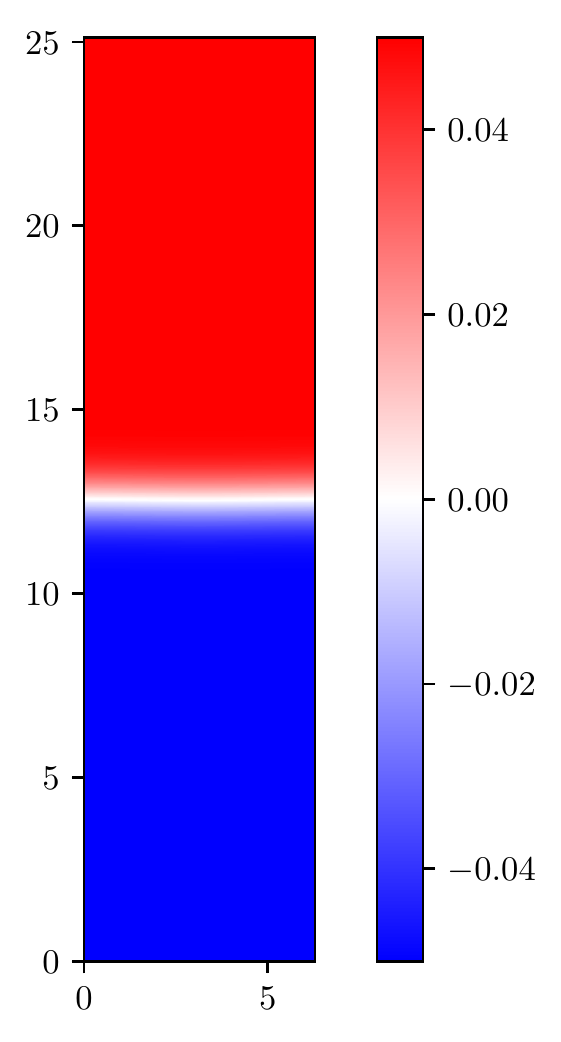}}
\subfloat[$T = 20.0$]{\includegraphics[width=0.2\textwidth]{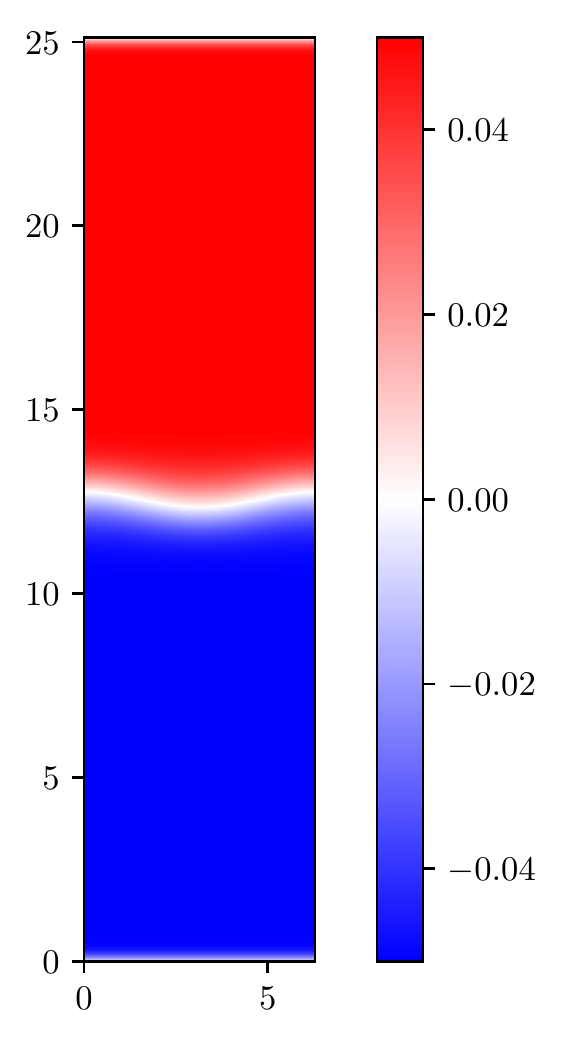}}
\subfloat[$T = 40.0$]{\includegraphics[width=0.2\textwidth]{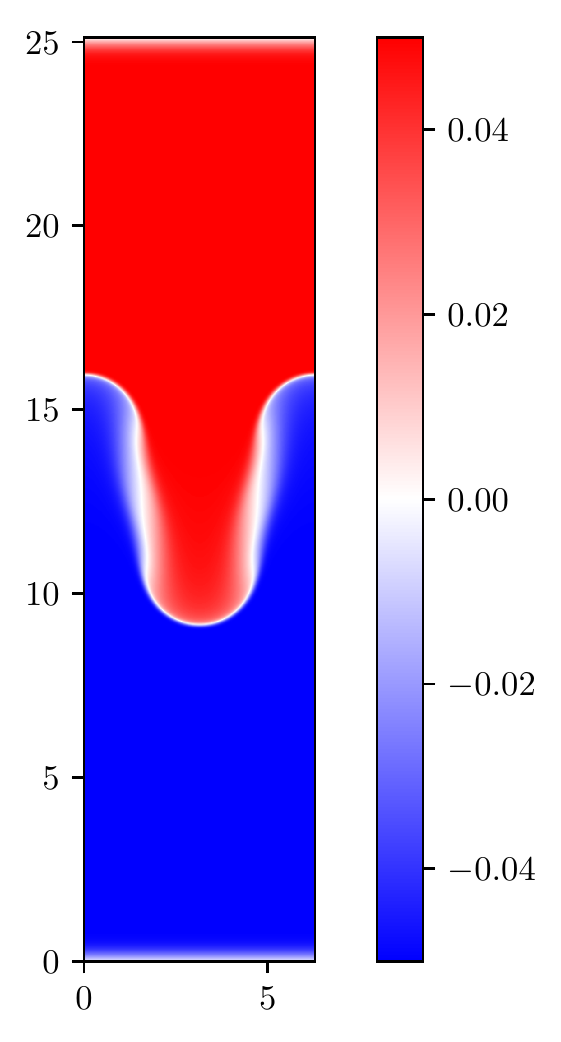}}
\subfloat[$T = 60.0$]{\includegraphics[width=0.2\textwidth]{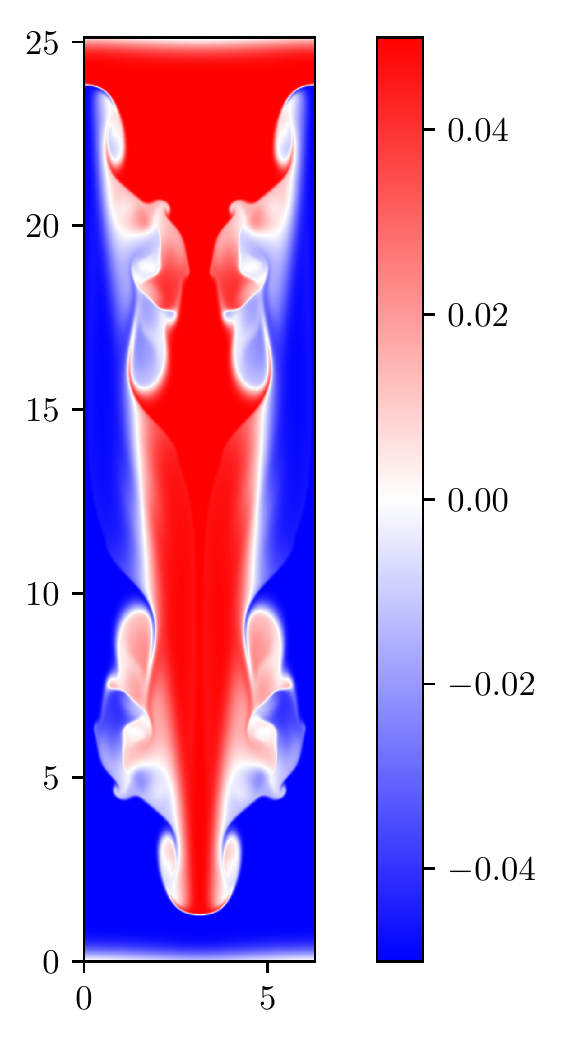}}
\caption{Slice of density-field of Rayleigh--Taylor instability on $256\times256\times 1024$ grid}
\label{fig:rti3ddensity}
\end{figure}

Here a $256\times256\times1024$ grid was used with $\nu = 10^{-3}$. These
parameters are similar to those used in \cite{young2001miscible}, though our
simulation differs in that we use a periodic domain: $\zeta(x,y) =
0.01\cos(x)\cos(y)$.
We used a reference timestep of $2\times 10^{-3}$.
As in the $2D$-case, we ran a reference simulation on a finer
$(512\times512\times2048)$ grid and compared the results to those on the
coarser grid. For both density and velocity, the (relative)
$L_2$-norm of the difference was less than $1\%$. In the maximum-norm, the
(relative) norm of the difference was $2.5\%$ and $20\%$ for the velocity field
and the density, respectively. Differences of the rate of dissipation of kinetic
energy and of the kinetic energy were both on the order of $10^{-10}$.
Differences in bubble-height were within the resolution of the grid. The
simulation was run from time $T=0$ to $T=60.0$. A plot of the density field
of the solution is plotted at various times in Figure \ref{fig:rti3ddensity}.

\begin{figure}
      \centering
      %\subfigure{\includegraphics[width=0.8\textwidth]{Figures/combined_plots_rt_3d/udiff_vs_f_evaluations_edited.pdf}}
      \includegraphics[width=0.5\textwidth]{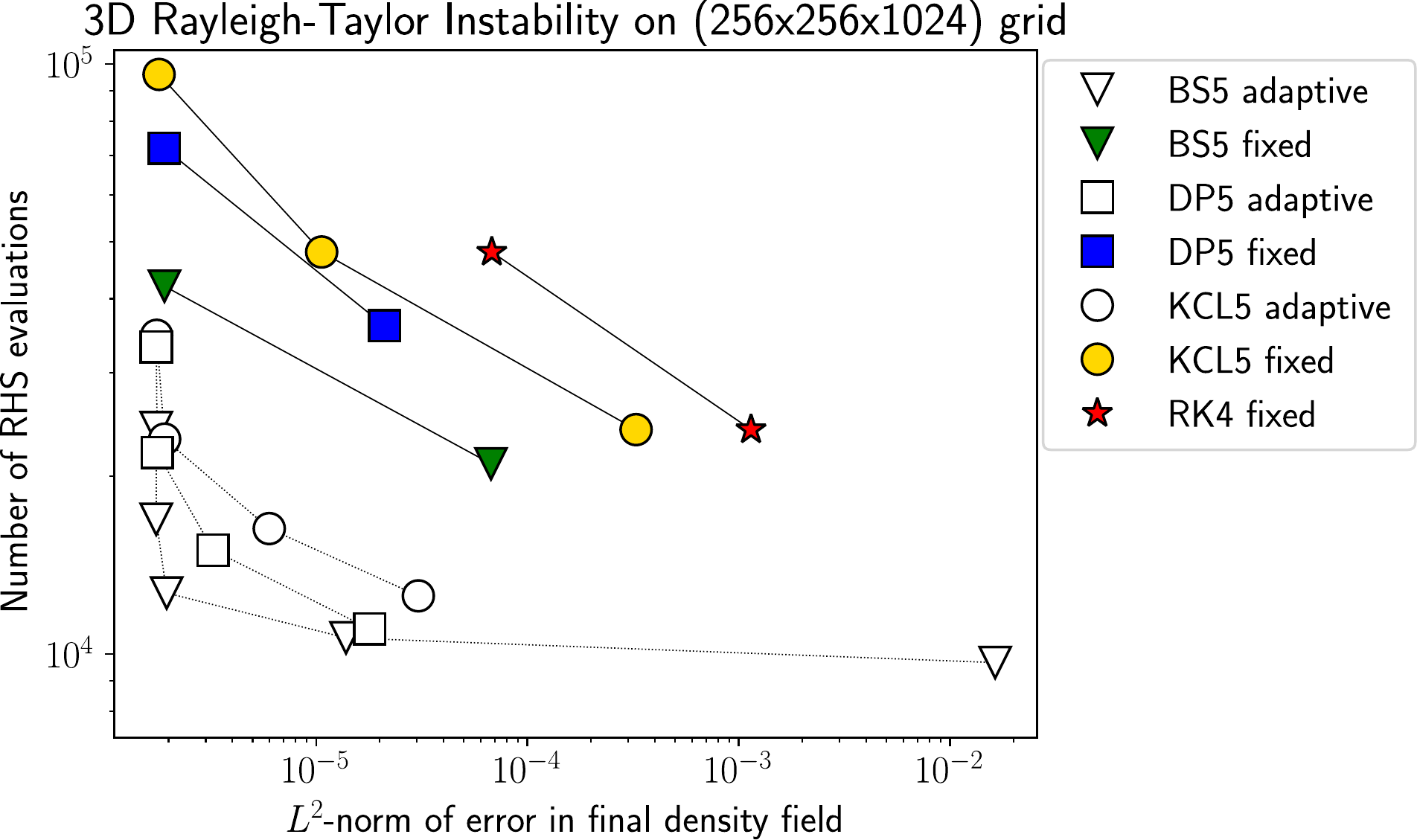}
      \includegraphics[width=0.5\textwidth]{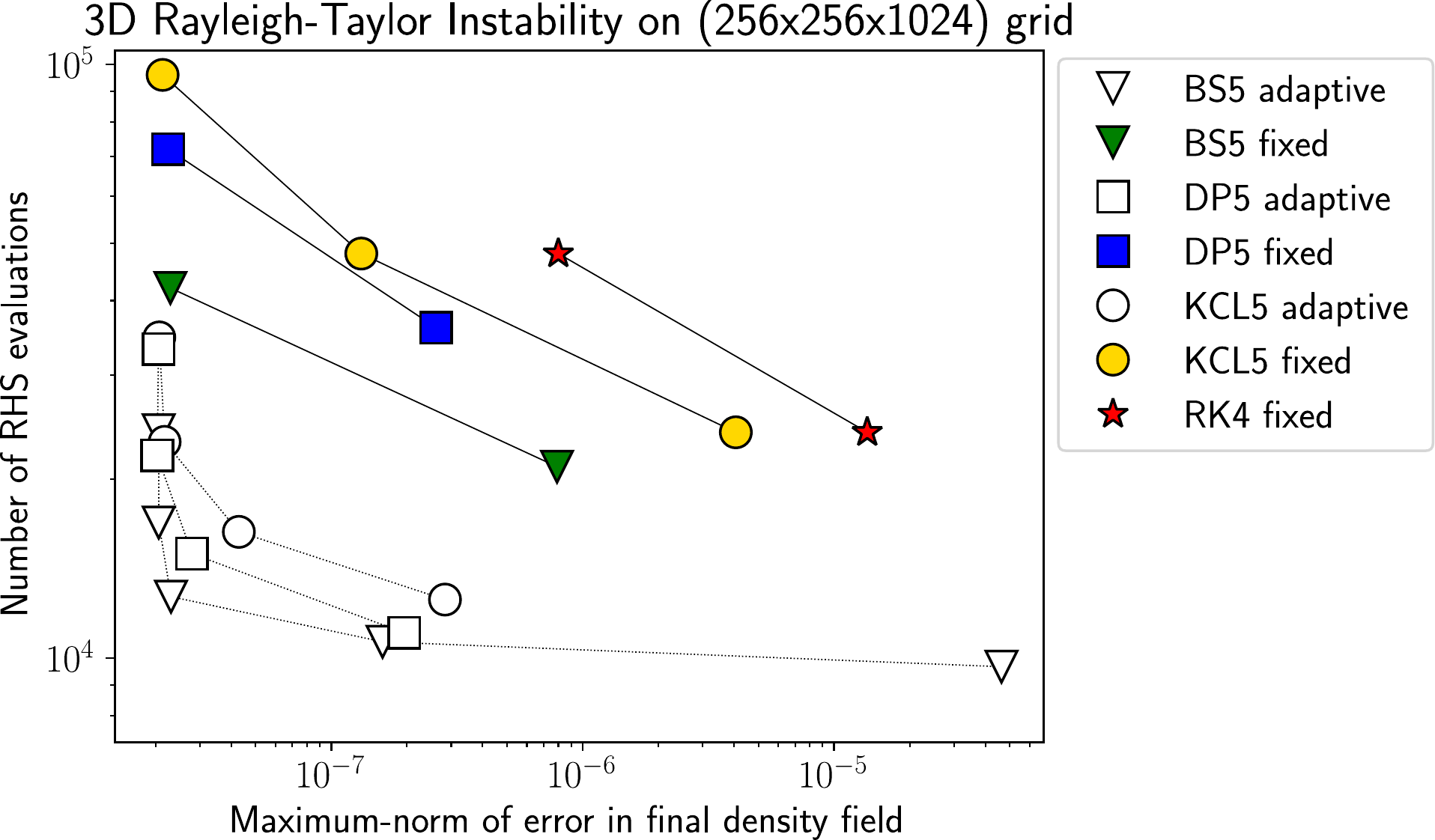}
      \caption{Comparison of time integration methods at $\nu^{-1}=1000$}\label{RTI3D}
\end{figure}

We used a range of tolerances from $10^{-8}$ to $10^{-4}$. Results for all methods are shown in Figure \ref{RTI3D}.
All the scheme with adaptive timestep outperform the schemes with fixed timestep, except for the
BS5 scheme with a fixed timestep whose number of RHS evaluations are close to those
of the KCL5.
Among the methods with adaptive timestep, the BS5 scheme is the most efficient for this problem.

\subsection{Homogeneous isotropic turbulence}

%\begin{table}
%\begin{tabular}{cccccc}
% \hline
%Paper/code & Machine & Grid & Largest Re & Time integrator & Step size \\ \hline
%Yokokawa et. al. 2002 \cite{yokokawa200216} & Earth Simulator & $4096^3$ & & RK4 & Not specified \\
%Ishihara et. al. 2007 \cite{ishihara2007small} & Earth Simulator & $4096^3$ & 36,500 & RK4 & $0.25e-3$ to $1e-3$ \\
%Donzis et. al. 2008 \cite{donzis2008dissipation} & Various & $2048^3$ & & RK2 w/int. factor & CFL 0.5 \\
%Yeung et. al. 2012 \cite{yeung2012dissipation} & TACC & $4096^3$ & $1000^*$ & RK2 w/int. factor & CFL 0.6 \\
%Yeung et. al. 2015 \cite{yeung2015extreme} & Blue waters & $8192^3$ & 45,000 & RK2 w/int. factor & CFL 0.6 \\
%%Tarang & & & & Euler, RK2, RK4 & \\
%%Turbo & & & & RK3 & \\
%%Philofluid & & & & 4th-order low-storage RK & \\
%%hit3d & & & & Euler, AB2 & \\
%\hline
%\end{tabular}
%\caption{Time stepping methods used in Fourier pseudo-spectral DNS of homogeneous isotropic turbulence.\label{history}}
%\end{table}

%The test-cases so far have all had some kind of symmetry in the initial conditions.
%\comment{Probably should say something about Homogeneous isotropic Turbulence here}.
Generally speaking, turbulent flows consist of vortices of various scales interacting with each other.
Energy is transferred from vortices of larger scale to
energy-dissipative vortices of smaller scale. This physical process is known as the turbulent energy cascade.
In this context, the study of homogeneous, isotropic turbulence (HIT) is very important for two reasons:
i) the smallest turbulent structures in most
turbulent flows have an almost isotropic behaviour and therefore it is hoped that these small structures,
often not represented in numerical simulations, can be modelled correctly, ii) it is possible to study and
understand an important part of HIT analytically.
We are particularly interested in proposing improved time integration for HIT since some of the largest
DNS runs ever have been devoted to this problem 
\cite{yokokawa200216,ishihara2007small,donzis2008dissipation,yeung2012dissipation,yeung2015extreme}.

We simulated homogeneous isotropic turbulence on the cube $[0,2\pi]^3$. The initial conditions used were based on \cite{sullivan94} and were initialized pointwise with each component at each point having random phase and absolute value (i.e., amplitude) given by
\begin{align}
    \frac{C|k|2 N^3}{(2\pi)^3}\exp\left(\frac{-|k|^2}{a^2}\right),
\end{align}
where $k$ is the wavenumber, and $a = 9.5$. The constant $C$ was chosen so that the total kinetic energy was of a similar magnitude as in the Taylor--Green vortex for the $256^3$ and $512^3$ case. We subtracted $\frac{(k\cdot \hat{u})k}{|k|^2}$ to make the initial conditions divergence-free.
\subsubsection{Forcing}
In order to replace the energy lost, we used a method of determinstic forcing
after each timestep based on the approach presented in \cite{sullivan94}. Here the lower
wavenumbers (with $0 < |k| \leq k_f$) are scaled by an appropriate factor to
compensate for energy lost by the system to diffusion.

\begin{figure}
      \centering
      \includegraphics[width=0.5\textwidth]{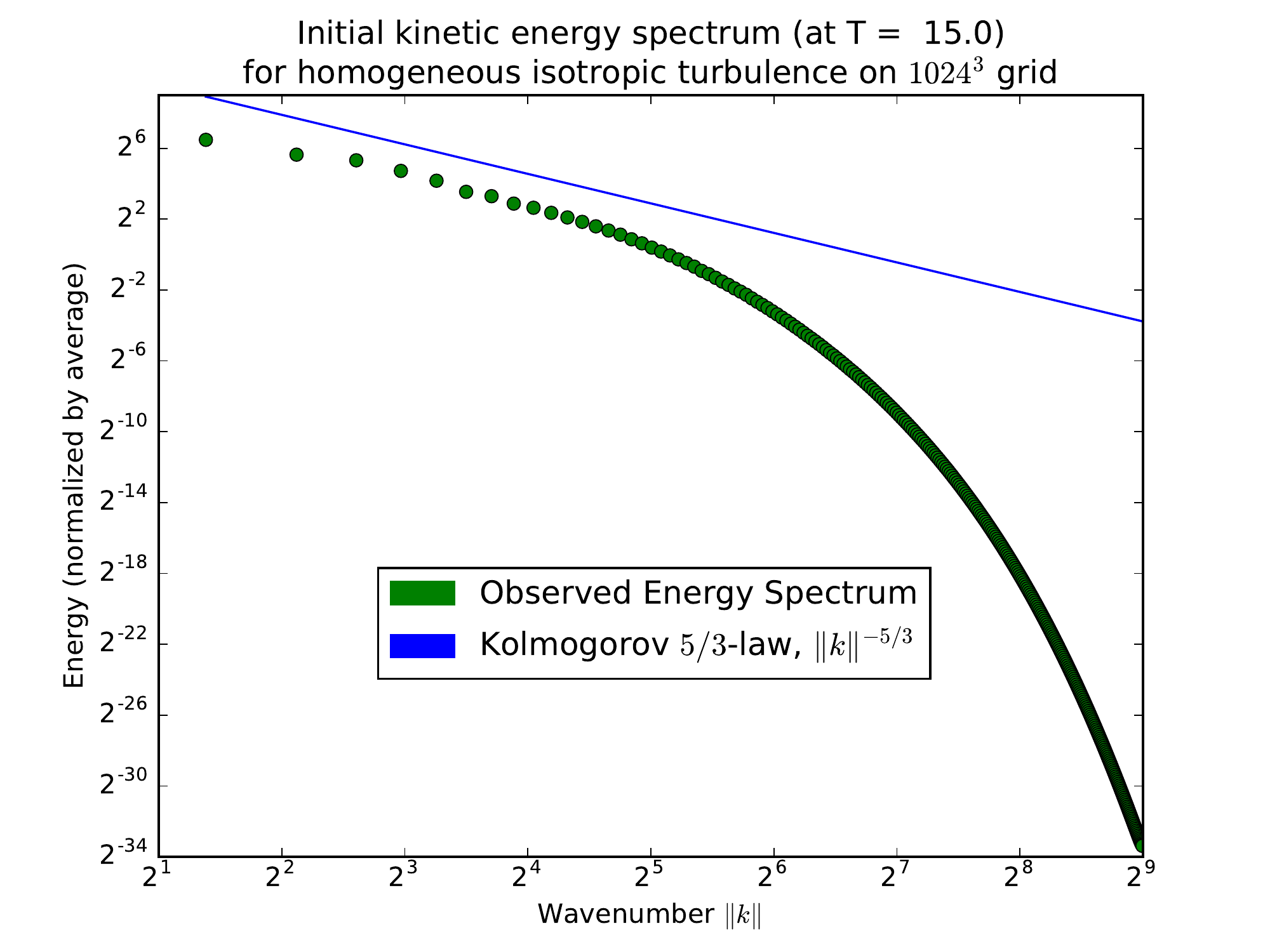}
      \caption{Initial energy spectrum for $\nu^{-1}=2000$ HIT run}\label{HIT1024spec}
      \label{fig:energyspectra}
\end{figure}

We calculated the energy spectra of the initial velocity fields like done in \cite{lamorgese2005direct} . Figure \ref{fig:energyspectra} shows the energy spectrum for the $1024^3$ run, which had $\nu^{-1} = 2000$ and $k_f = 8$. 

To generate an appropriate initial condition, the simulation was run from $T=0$
to $T=15.0$ using the BS5 time-integrator with adaptive stepsize (with
tolerance $10^{-7}$ on the $256^3$ run and $10^{-6}$ on the others).
The resulting energy spectrum is shown in Figure \ref{HIT1024spec}.  We see that
the high-order adaptive integrator yields the expected statistical state. 
After reaching this state, the adaptive timestepping gave a step
size that corresponds to a CFL number of about $0.6$ for the $1024^3$ run.

\subsubsection{Decaying turbulence test runs}

\begin{figure}
\centering
\subfloat[$256^3$ grid]{\includegraphics[width=0.3\textwidth]{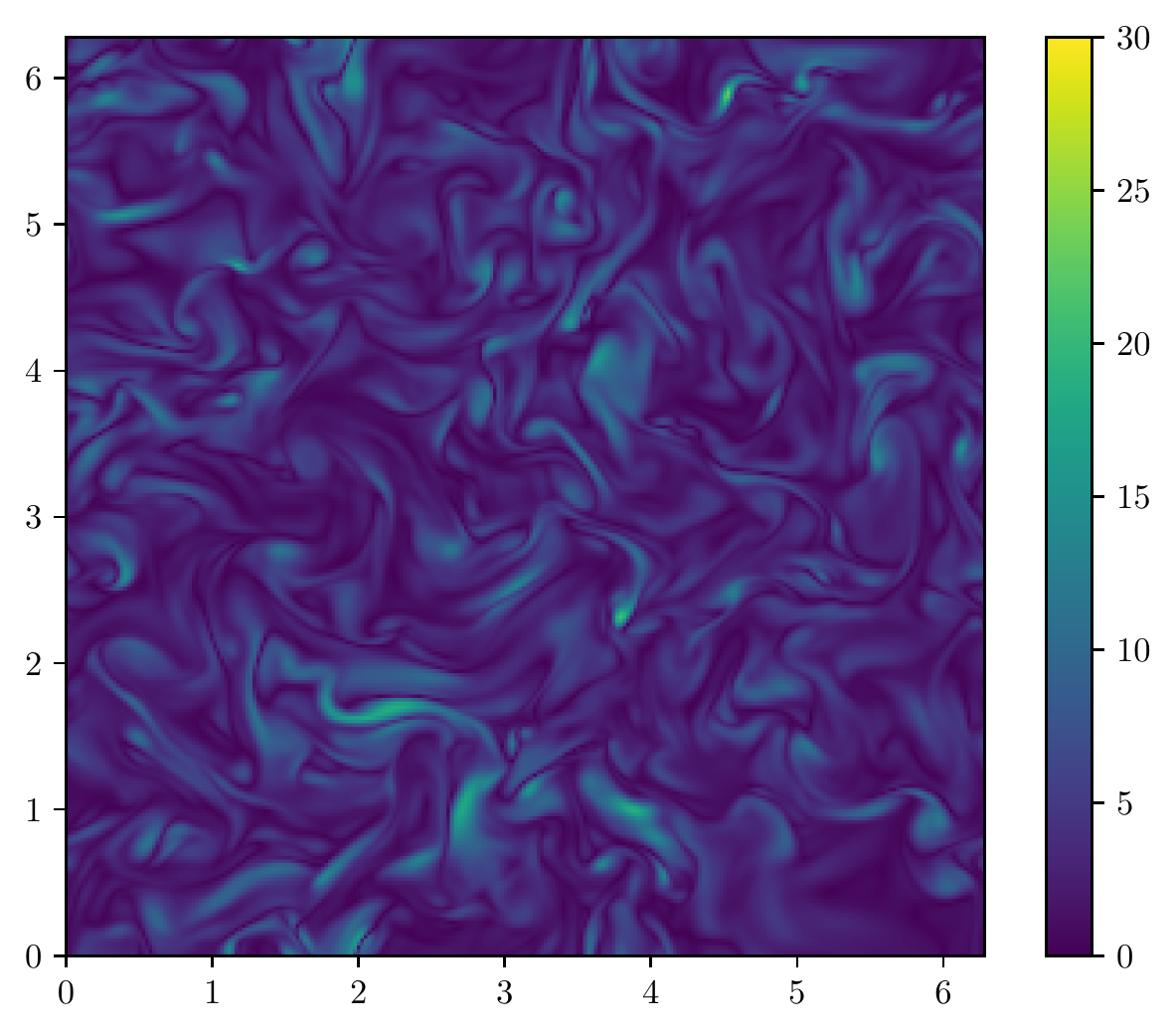}}
\subfloat[$512^3$ grid]{\includegraphics[width=0.3\textwidth]{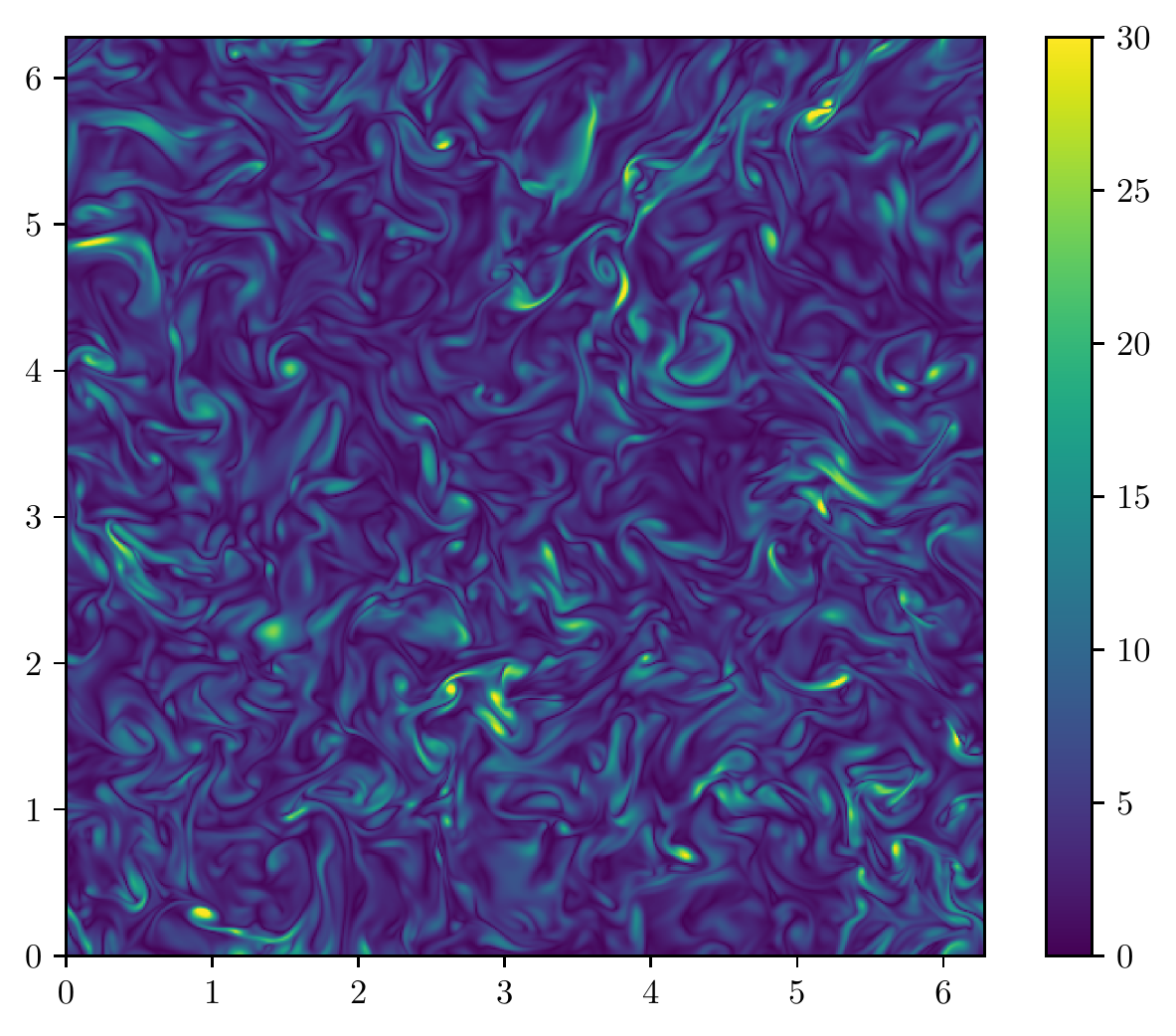}}
\subfloat[$1024^3$ grid]{\includegraphics[width=0.3\textwidth]{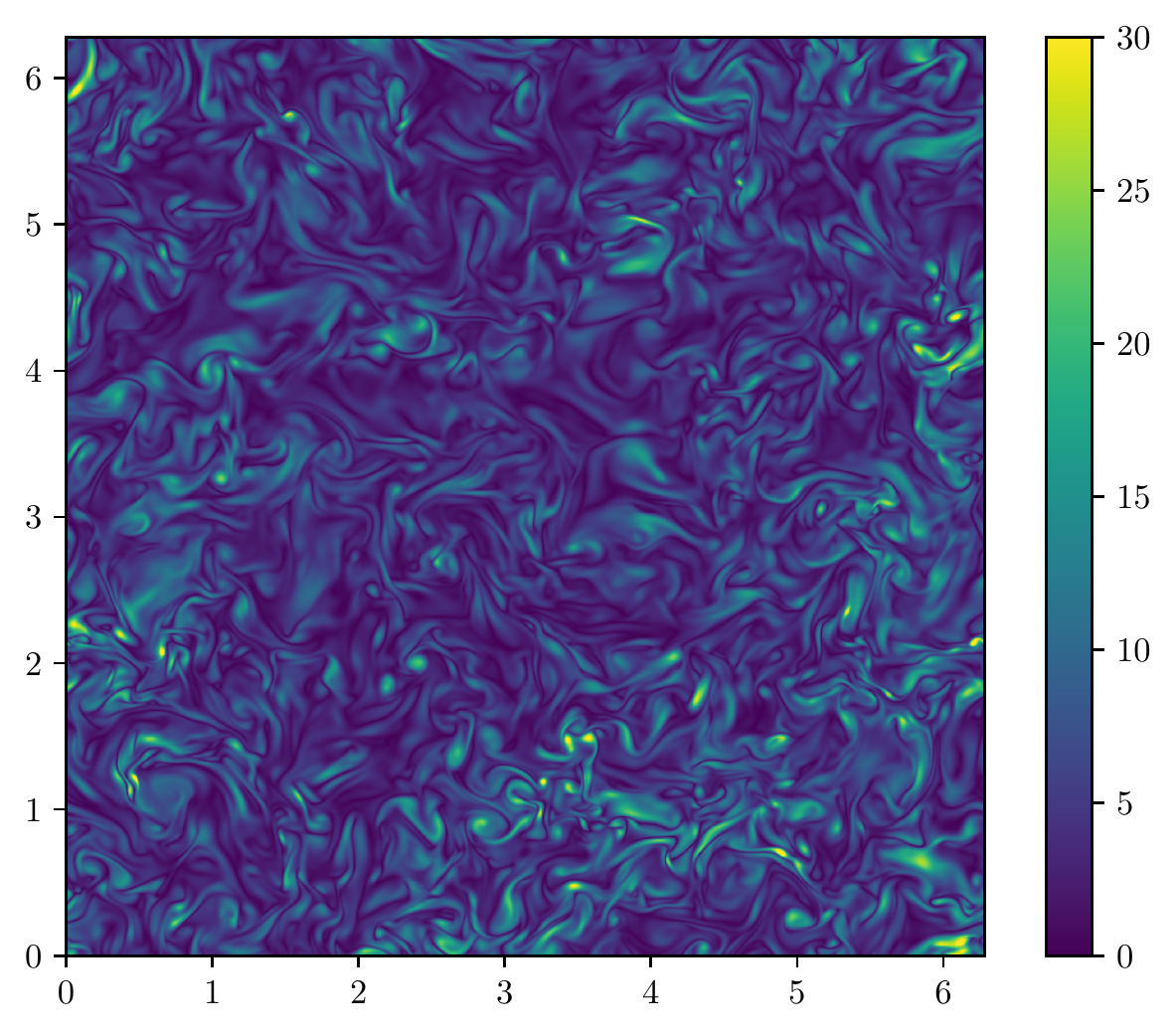}}
\caption{Vorticity slice ($x = 0$) at $t=15$}
\label{fig:hitslice}
\end{figure}

To test each integrator, the simulation was restarted from $t=15$ and
run for a short time without forcing. A slice of the vorticity of the solution at $t = 15$ is shown in Figure \ref{fig:hitslice}.

We ran the simulation on various grid sizes with $\nu$ chosen based on
the value of $\nu$ for the Taylor--Green vortex and on
\cite{rosa2016settling}.
As a very stringent test of the time integrators, we compare pointwise
solution values of the resulting turbulent flow field.\footnote{
After running these (expensive) simulations it was found that the zero and
Nyquist frequencies of the first Fourier transform had a non-zero imaginary
part. This could
slightly reduce the step size chosen by adaptive time-stepping, in a uniform way
for all methods. The energy spectrum in Figure \ref{fig:energyspectra} includes only the real part.  }

\begin{figure}
      \centering
      \includegraphics[width=0.8\textwidth]{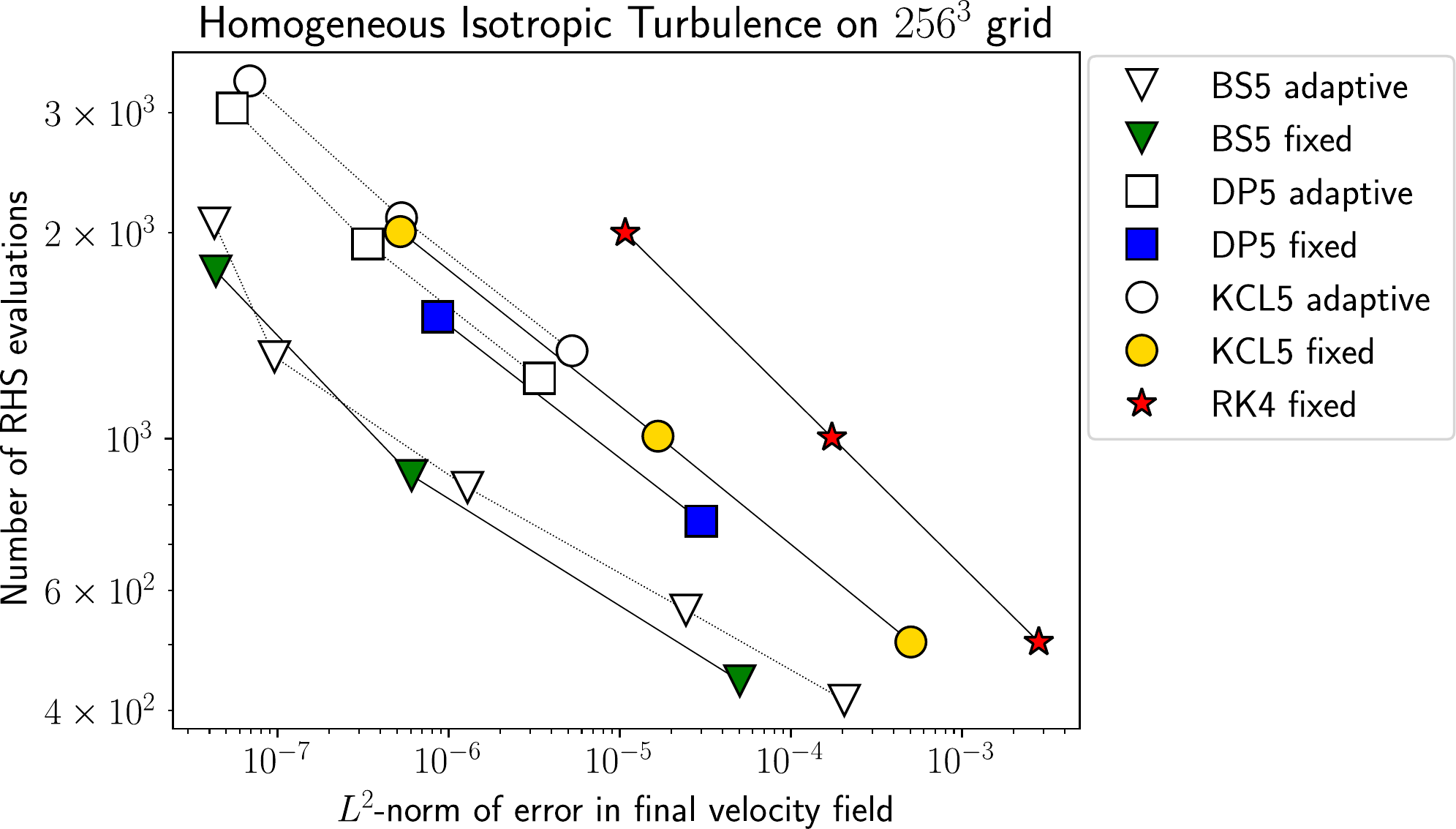}
      \includegraphics[width=0.8\textwidth]{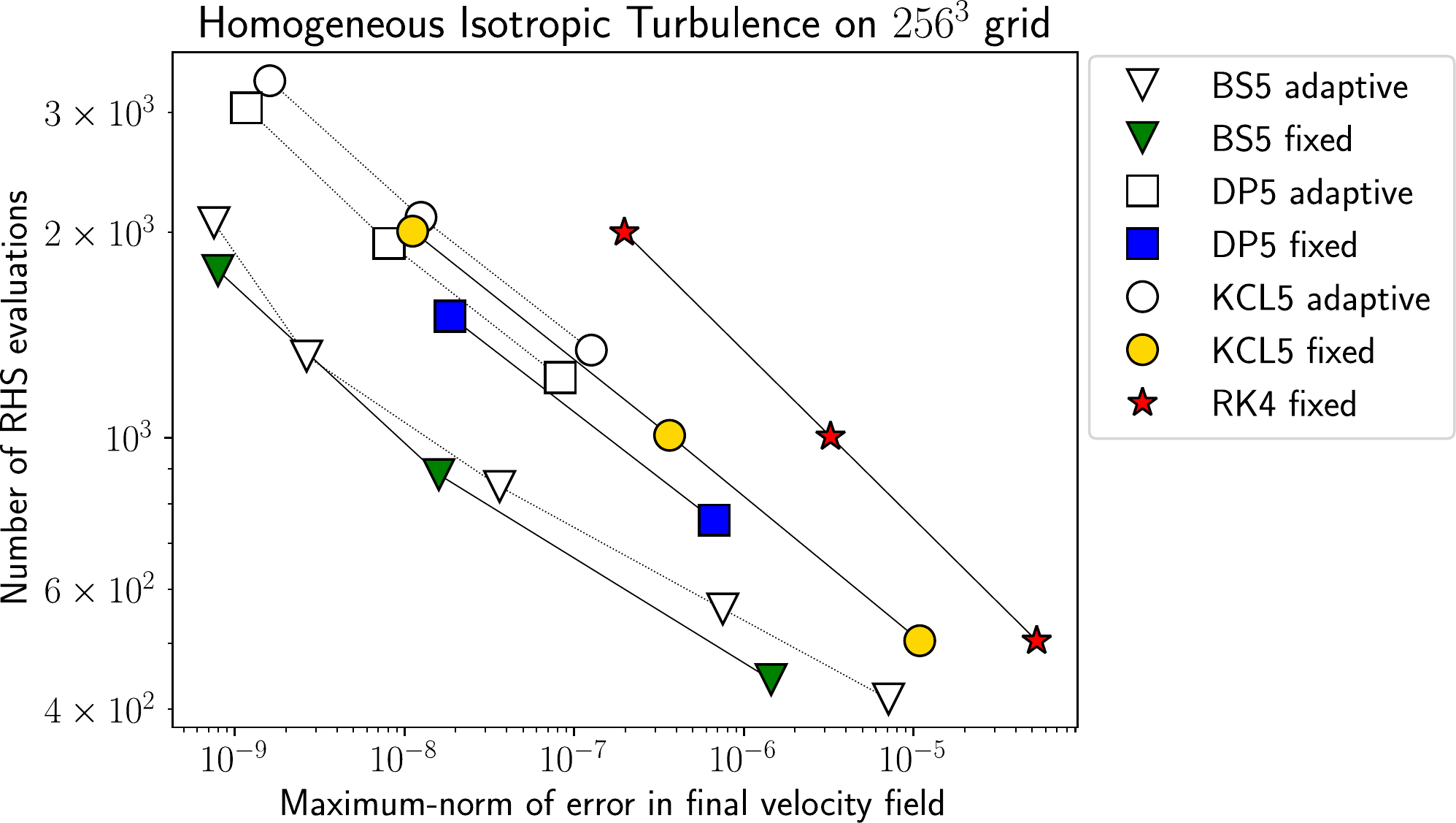}
      \caption{Comparison of time integration methods at $\nu^{-1}=800$}\label{HIT256}
\end{figure}

Results for a $256^3$ grid are shown in Figure \ref{HIT256}. Here we used $\nu^{-1} = 800$ and
$k_f = 4\sqrt{2}$. The simulation was run from $T=15.0$ to $T=16.0$. As a
reference we used RK4 with timestep $10^{-3}$. We used a range of tolerances between $10^{-8}$ and $10^{-4}$.

\begin{figure}
      \centering
      \includegraphics[width=0.8\textwidth]{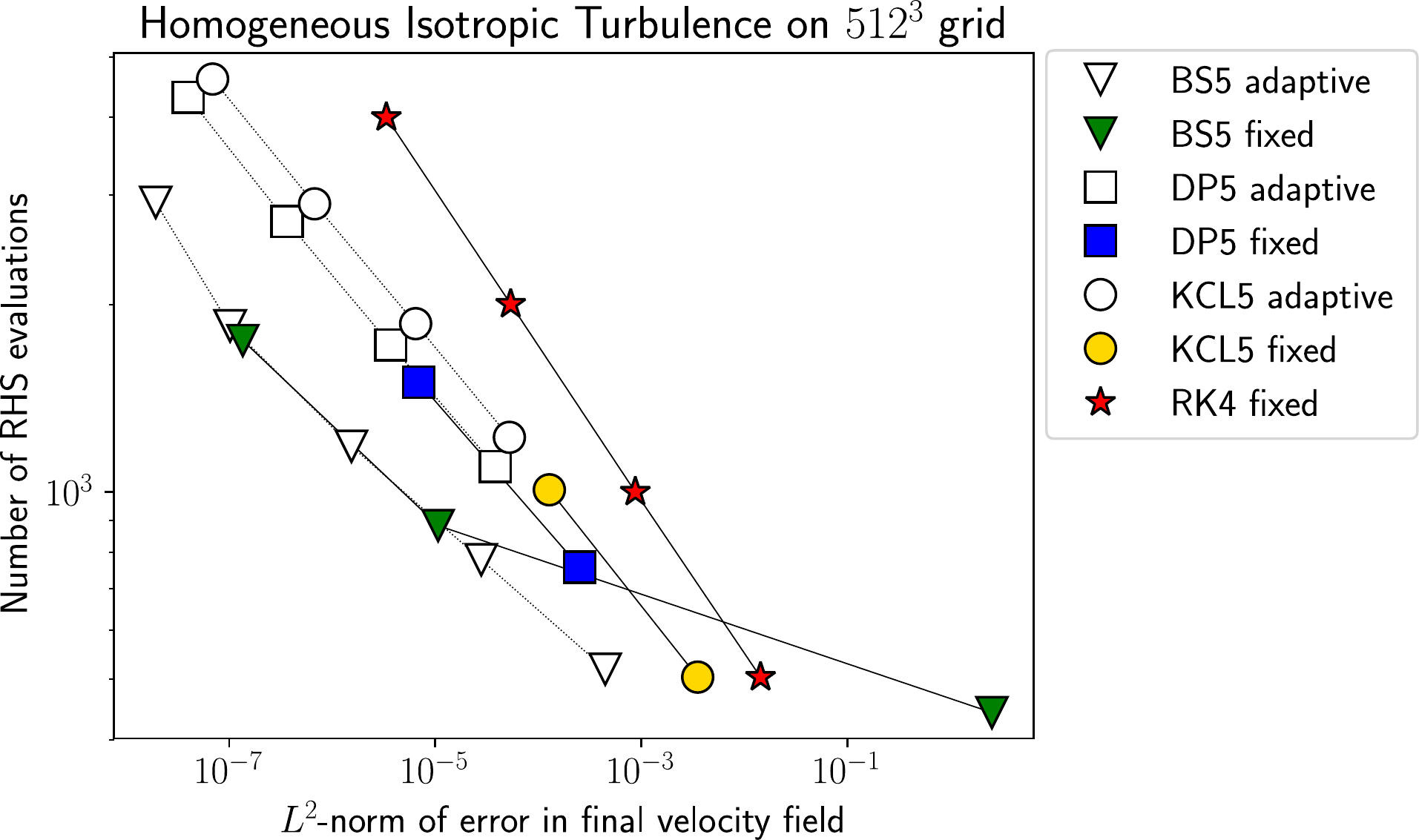}
      \includegraphics[width=0.8\textwidth]{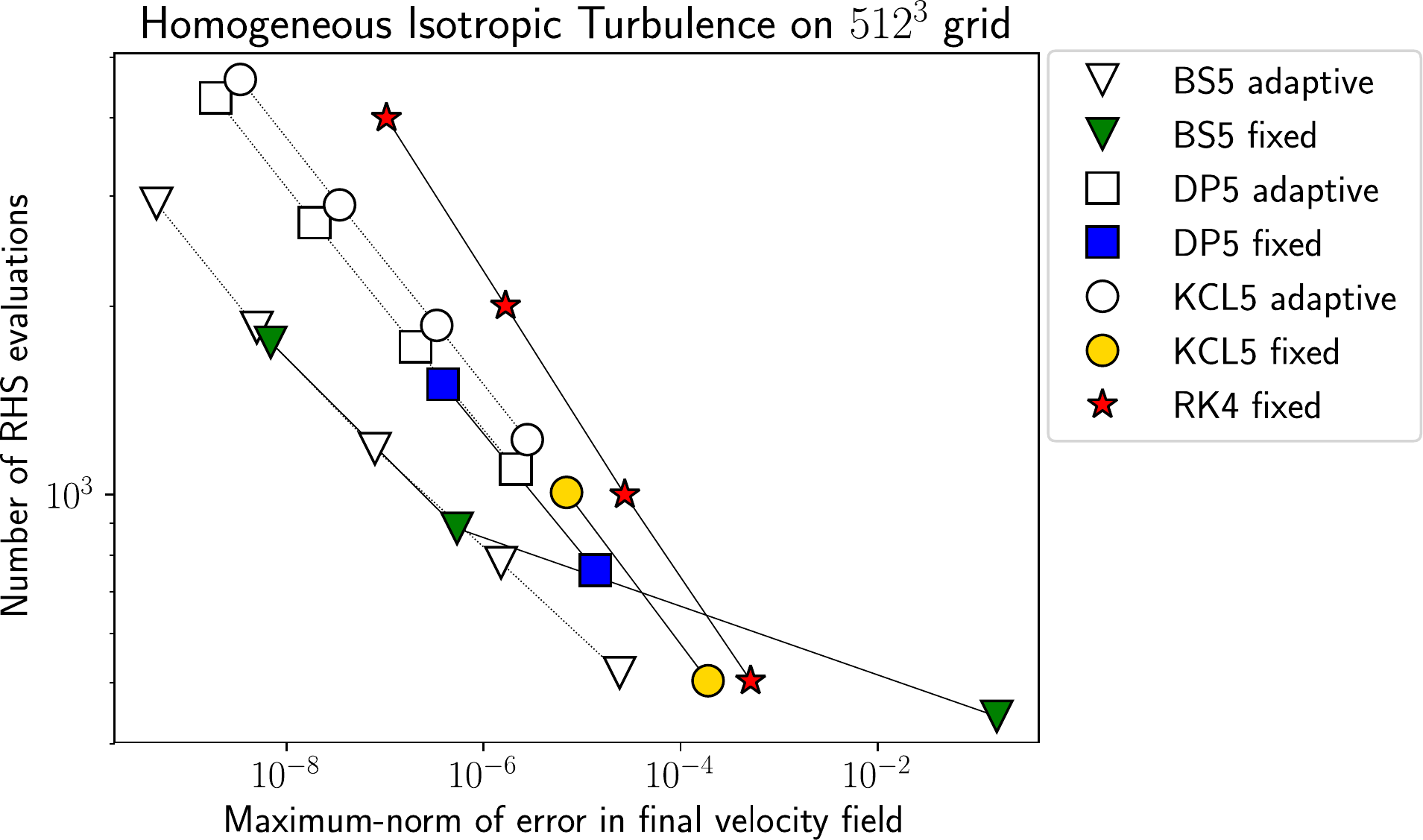}
      \caption{Comparison of time integration methods at $\nu^{-1}=1600$}\label{HIT512}
\end{figure}

Results for a $512^3$ grid are shown in Figure \ref{HIT512}.  Here we used $\nu^{-1} = 1600$ and
$k_f = 8$. The simulation was run from $T=15.0$ to $T=15.5$. As a reference we
used RK4 with timestep $5\times10^{-4}$. We used a range of tolerances between $10^{-8}$ and $10^{-4}$.

\begin{figure}
      \centering
      \includegraphics[width=0.8\textwidth]{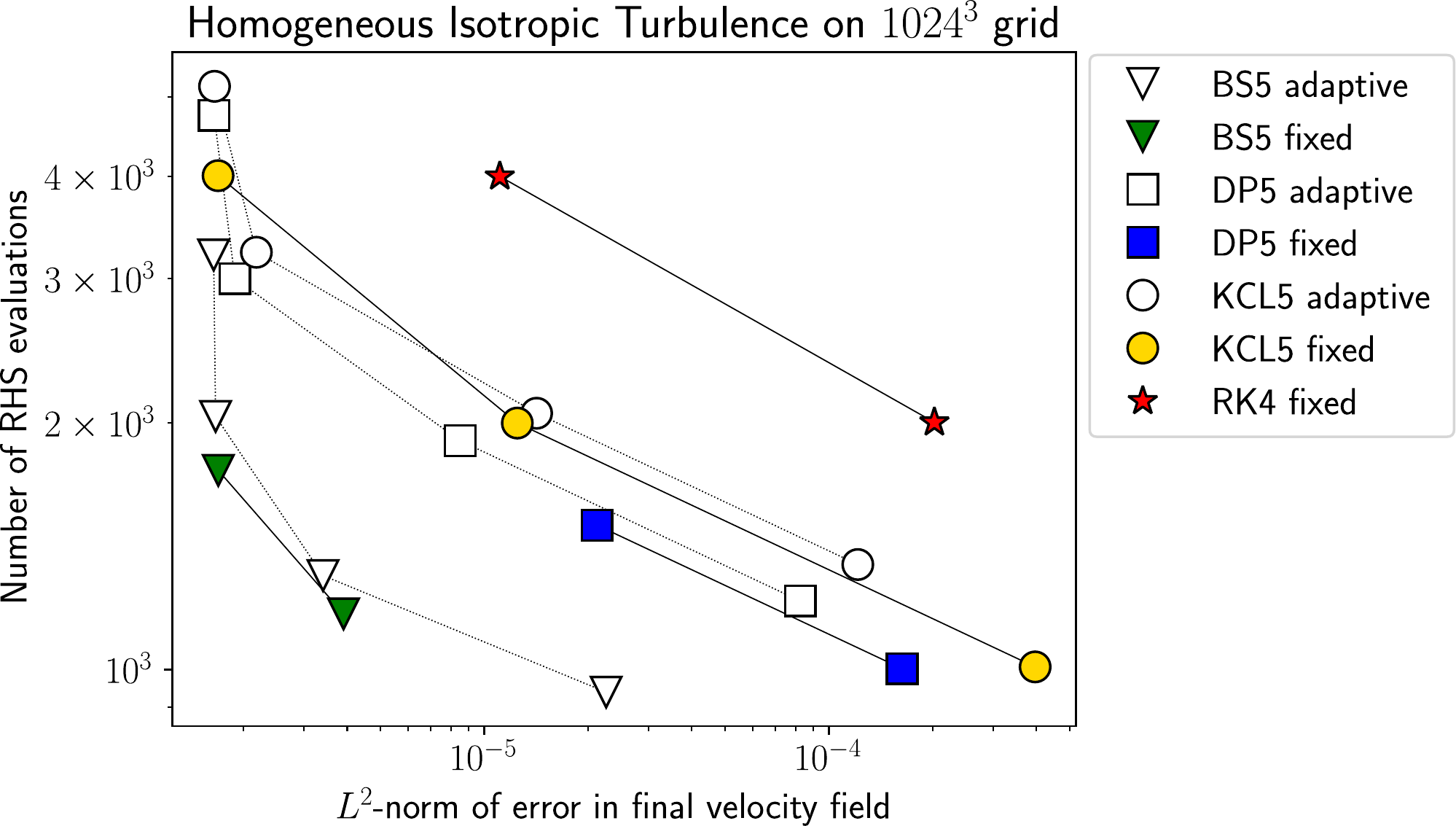}
      \includegraphics[width=0.8\textwidth]{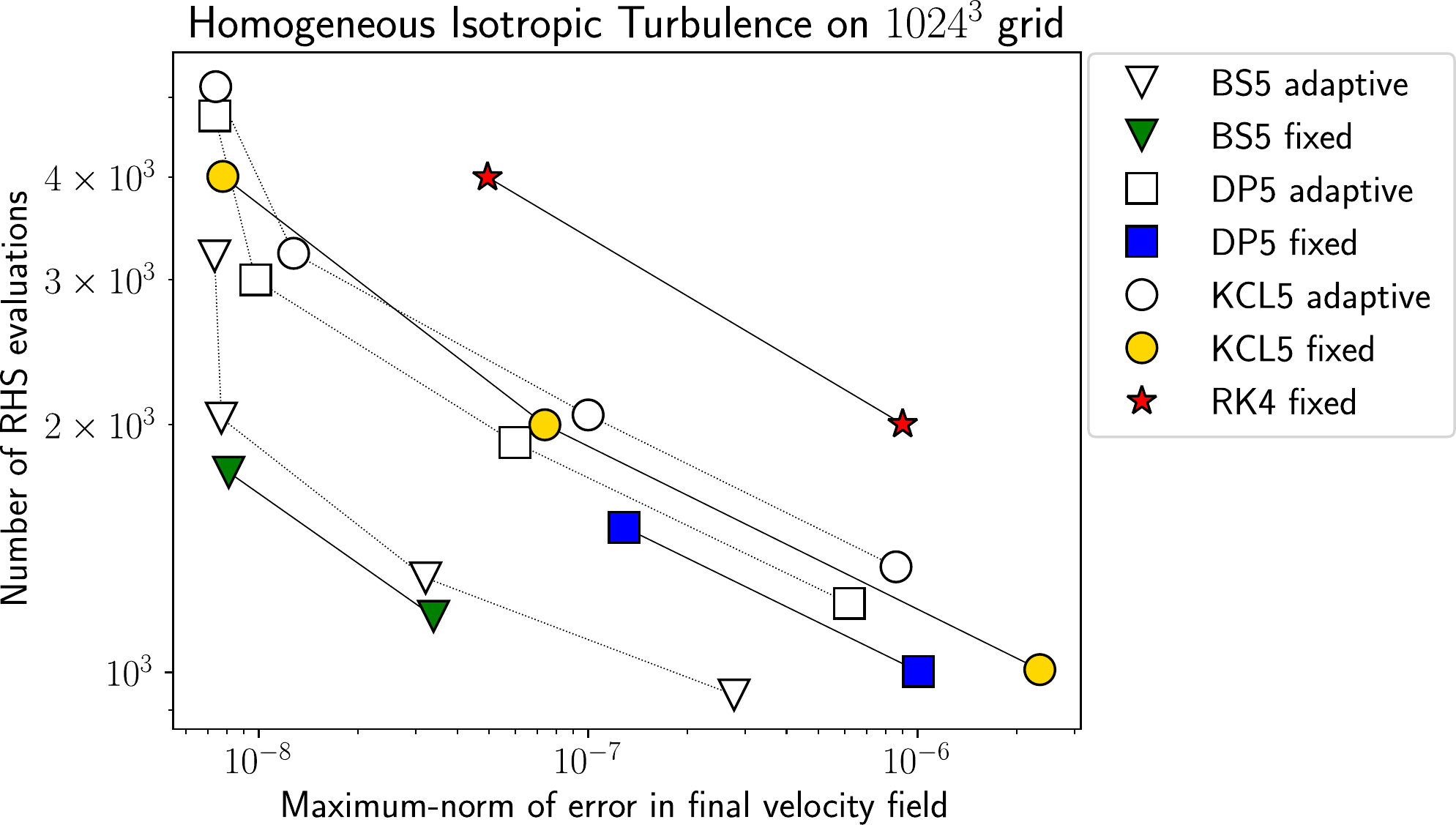}
      \caption{Comparison of time integration methods at $\nu^{-1}=2000$\label{HIT1024}}
\end{figure}

Results for a $1024^3$ grid are shown in Figure \ref{HIT1024}.  Here we used $\nu^{-1} = 2000$ and $k_f = 8$. The simulation was run from $T=15.0$ to $T=15.5$. As a reference we used RK4 with timestep $5\times10^{-4}$.
We used a range of tolerances between $10^{-8}$ and $10^{-5}$.

As for the other test cases presented herein, the BS5 method with adaptive timestep requires less RHS function evaluations. However, for this flow problem, 
the BS5 scheme with a fixed timestep requires a number of RHS evaluations very
close if not better than the BS5 with automatic timestep selection.   

\section{Discussion}
While one should be cautious about drawing general results from a small selection
of test problems, our results consistently suggest that:
\begin{itemize}
    \item High-order, adaptive time stepping yields accurate results for pseudospectral DNS of turbulent flow.
            The use of significantly larger time step sizes with such methods does not negatively affect the resolution of
            fine-scale features of the flow when compared with the use of
            lower-order methods (and correspondingly smaller step sizes).
    \item Optimized fifth-order methods, especially the BS5 pair, are more efficient than lower-order
            methods whenever moderate to high accuracy is desired.  Specifically, similar accuracy can
            be obtained at a cost that is reduced by 2x to 10x.
    \item Error estimation and step size adaptivity can be highly beneficial for problems involving the development of
            instability from an initially laminar flow.  It can also be useful in automatically finding
            an appropriate step size, even if the characteristics of the flow do not change significantly
            during the simulation.
\end{itemize}

Similar work remains to be done for the development of time integrators for incompressible flow
in the presence of boundaries.
Many questions remain to be investigated regarding the application of modern time discretizations
to pseudo-spectral DNS.  In particular, we expect that exponential methods (wherein the linear
diffusive terms are handled directly via the matrix exponential) may yield even more substantial
improvements.  Improved step size controllers such as those of \cite{soderlind2006time} will also
likely show improvements over the more standard techniques used here.

\section{Acknowledgements}
This research used the resources of the Supercomputing Laboratory and Extreme Computing Research Center at the King Abdullah University of Science \& Technology (KAUST) in Thuwal, Saudi Arabia. N.S. was supported by the KAUST Visiting Student Research Program. N.S. also acknowledges support from the Priority Programme SPP1881 Turbulent Superstructures of the Deutsche Forschungsgemeinschaft.
M. M. acknowledges support from the 4DSpace Strategic Research Initiative at the University of Oslo. 
\bibliography{bibliography}

\appendix

\section{Coefficients of Runge--Kutta pairs}
For the convenience of the reader, we give here the coefficients of each Runge--Kutta pair
used in this work.

\subsection{Dormand--Prince 5(4) pair (DP5)}
{ % begin box to localize effect of arraystretch change
\renewcommand{\arraystretch}{1.5}
\begin{align}
\begin{array}{c|ccccccc}
 &  &  &  &  &  &  & \\
\frac{1}{5} & \frac{1}{5} &  &  &  &  &  & \\
\frac{3}{10} & \frac{3}{40} & \frac{9}{40} &  &  &  &  & \\
\frac{4}{5} & \frac{44}{45} & - \frac{56}{15} & \frac{32}{9} &  &  &  & \\
\frac{8}{9} & \frac{19372}{6561} & - \frac{25360}{2187} & \frac{64448}{6561} & - \frac{212}{729} &  &  & \\
1 & \frac{9017}{3168} & - \frac{355}{33} & \frac{46732}{5247} & \frac{49}{176} & - \frac{5103}{18656} &  & \\
1 & \frac{35}{384} &  & \frac{500}{1113} & \frac{125}{192} & - \frac{2187}{6784} & \frac{11}{84} & \\
\hline
 & \frac{35}{384} &  & \frac{500}{1113} & \frac{125}{192} & - \frac{2187}{6784} & \frac{11}{84} & \\
 & \frac{5179}{57600} &  & \frac{7571}{16695} & \frac{393}{640} & - \frac{92097}{339200} & \frac{187}{2100} & \frac{1}{40}
\end{array}
\end{align}
}
The first set of weights listed is used to advance the solution.  The second set ($\hat{b}$) is used
only for error estimation.

\subsection{KCL5 low-storage scheme}
Because of the special structure and desirable low-storage implementation of this scheme,
the Butcher tableau is not the most useful format for its presentation.  Instead, we refer the
reader to \cite[p.~190]{Kennedy2000} where the low-storage coefficients are presented.

\subsection{Bogack--Shampine 5(4) pair (BS5)}
{ % begin box to localize effect of arraystretch change
\renewcommand{\arraystretch}{1.5}
\begin{align}
\begin{array}{c|cccccccc}
 &  &  &  &  &  &  &  & \\
\frac{1}{6} & \frac{1}{6} &  &  &  &  &  &  & \\
\frac{2}{9} & \frac{2}{27} & \frac{4}{27} &  &  &  &  &  & \\
\frac{3}{7} & \frac{183}{1372} & - \frac{162}{343} & \frac{1053}{1372} &  &  &  &  & \\
\frac{2}{3} & \frac{68}{297} & - \frac{4}{11} & \frac{42}{143} & \frac{1960}{3861} &  &  &  & \\
\frac{3}{4} & \frac{597}{22528} & \frac{81}{352} & \frac{63099}{585728} & \frac{58653}{366080} & \frac{4617}{20480} &  &  & \\
1 & \frac{174197}{959244} & - \frac{30942}{79937} & \frac{8152137}{19744439} & \frac{666106}{1039181} & - \frac{29421}{29068} & \frac{482048}{414219} &  & \\
1 & \frac{587}{8064} &  & \frac{4440339}{15491840} & \frac{24353}{124800} & \frac{387}{44800} & \frac{2152}{5985} & \frac{7267}{94080} & \\
\hline
 & \frac{587}{8064} &  & \frac{4440339}{15491840} & \frac{24353}{124800} & \frac{387}{44800} & \frac{2152}{5985} & \frac{7267}{94080} & \\
 & \frac{2479}{34992} &  & \frac{123}{416} & \frac{612941}{3411720} & \frac{43}{1440} & \frac{2272}{6561} & \frac{79937}{1113912} & \frac{3293}{556956}
\end{array}
\end{align}
}
The first set of weights listed is used to advance the solution.  The second set ($\hat{b}$) is used
only for error estimation.

\end{document}